\documentclass[11pt,reqno]{amsart}

\usepackage{amsfonts}
\usepackage{eurosym}
\usepackage{amssymb}
\usepackage{amsthm}
\usepackage{amsmath}
 \usepackage[foot]{amsaddr}
\usepackage{bm}
\usepackage{cite}
\usepackage{mathrsfs}
\usepackage[usenames,dvipsnames]{xcolor}
\usepackage[utf8]{inputenc}
\usepackage[left=2cm, right=2cm, top=3cm]{geometry}
\usepackage{hyperref}
\hypersetup{colorlinks=true, linkcolor=blue, citecolor=red}

\usepackage{times}

\usepackage[shortlabels]{enumitem}
\usepackage{enumitem}
\newcommand{\subscript}[2]{$#1 _ #2$}

\numberwithin{equation}{section}
\usepackage[square,numbers,sectionbib]{natbib}
\usepackage{varioref}
\usepackage{enumitem}
\usepackage{listings}
\usepackage{color}
\definecolor{dkgreen}{rgb}{0,0.6,0}
\definecolor{gray}{rgb}{0.5,0.5,0.5}
\definecolor{mauve}{rgb}{0.58,0,0.82}

\allowdisplaybreaks
\newtheorem{theorem}{Theorem}[section]

\newtheorem{corollary}[theorem]{Corollary}
\newtheorem{definition}[theorem]{Definition}
\newtheorem{lemma}[theorem]{Lemma}
\newtheorem{proposition}[theorem]{Proposition}
\theoremstyle{remark}
\newtheorem{remark}[theorem]{Remark}

\def\d{{\rm d}}
\def \l {\langle}
\def \r {\rangle}

\def\uu{\textbf{\textit{u}}}

\def\E{\mathcal E}

\newcommand{\Fsec}{F^{\prime\prime}}
\newcommand{\numberset}{\mathbb}
\newcommand{\N}{\numberset{N}}
\newcommand{\R}{\numberset{R}}
\def\ry #1{{#1}}
\def\rx #1{{#1}}
\def\extra #1{{#1}}

\everymath{\displaystyle}

\def \au {\rm}
\def \ti {\it}
\def \jou {\rm}
\def \bk {\it}
\def \no#1#2#3 {{\bf #1} (#3), #2.}
\def \eds#1#2#3 {#1, #2, #3.}

\begin{document}

\title[The 3D strict separation property for the nonlocal Cahn-Hilliard equation with singular potential]
{The 3D strict separation property \\ for the nonlocal Cahn-Hilliard equation \\ with singular potential}
\author[The 3D strict separation property for the nonlocal Cahn-Hilliard equation with singular potential]{Andrea Poiatti$^*$}
\thanks{THE PRESENT RESEARCH HAS BEEN SUPPORTED BY MUR GRANT DIPARTIMENTO DI
ECCELLENZA 2023-2027.}
\address{$^*$Dipartimento di Matematica\\
Politecnico di Milano\\
Milano 20133, Italy}
\email{ andrea.poiatti@polimi.it}
\date{07 August 2022}
\keywords{Three-dimensional nonlocal Cahn-Hilliard equation, singular potential, strict separation property, regularization of weak solutions, convergence to equilibrium.}

\subjclass[2020]{35B40, 35B65, 35Q82, 35R09.}
\maketitle

\begin{abstract}
We consider the nonlocal Cahn-Hilliard equation with singular (logarithmic) potential and constant mobility in three-dimensional bounded domains and we establish the validity of the instantaneous strict separation property. This means that any weak solution, which is not a pure phase initially, stays uniformly away from the pure phases $\pm1$ from any positive time on. This work extends the result in dimension two for the same equation and gives a positive answer to the long standing open problem of the validity of the strict separation property in dimensions higher than two. In conclusion, we show how this property plays an essential role to achieve higher-order regularity for the solutions and to prove that any weak solution converges to a single equilibrium.
\end{abstract}

\section{Introduction}

The Diffuse Interface theory, also called Phase Field method, is one of the oldest and efficient approach to multi-phase problems. This approach is characterized by the notion of diffuse interface, meaning that the transition layer between the two phases or components has a narrow finite size. The interface is not explicitly tracked as in boundary integral and front-tracking methods. On the other hand, the phase state is incorporated into the macroscopic equations and the internal microstructures arise from the competition between the diffusion and aggregation mechanisms included in the free energy. The fundamental advantage of this theory is the natural 
representation of singular interfacial behaviors, such as topological change,
self-intersection, merger and pinch-off.

Consider a mixture of two incompatible substances A and B, which is homogeneously distributed and isothermal. Under certain circumstances, namely if the temperature
is above a critical threshold $\theta_c$, this configuration is stable; however, if suddenly
cooled down and kept at $\bar{\theta} < \theta_c$ , the initially (macroscopically) homogeneous alloy
evolves in a way such that A-rich and B-rich regions appear and grow. The Cahn-Hilliard equation was introduced in \cite{num5} and \cite{num13} to model this phenomenon in iron alloys, and it has now become a widespread model, since phase separation has become a paradigm also in Cell Biology (see, e.g., \cite{Dolgin}). Let $\Omega$ be a bounded domain in $\R^d$, $d=2,3$, filled with a binary solution consisting of A and B atoms, and let us fix a time horizon $T>0$. We define their relative mass fraction difference as $\phi$, which is the phase-field variable, whose  smooth but highly localized variation is associated with the (diffuse) interface. If the mixture is isothermal and the molar volume is uniform and independent on pressure, the system evolves in order to minimize the free energy functional
\begin{equation}
\mathcal{U}(\phi):=\int_\Omega \left(\frac{\epsilon}{2}|\nabla \phi|^2+\Psi(\phi)\right)dx,
\label{ener_loc}
\end{equation} 
where $\Psi(\phi)$ is the Helmholtz free energy density 
\begin{equation}
\Psi(s)=\frac{\bar{\alpha}}{2}((1+s)\text{ln}(1+s)+(1-s)\text{ln}(1-s))-\frac{\alpha_0}{2}s^2=F(s)-\frac{\alpha_0}{2}s^2,\ \ \ \forall s\in[-1,1],
\label{potential_loc}
\end{equation}
with $\bar{\alpha}$ such that $0<\bar{\alpha}<\alpha_0$, constants related to the temperature of the mixture.  The term $\epsilon$ is called capillary coefficient, related to the thickness of interfaces. The potential defined in this way is called \textit{singular}, whereas many authors (see, e.g., \cite{Gerbi14}) considered a proper approximation, which avoids the fact that $\Psi'$ is unbounded at the pure phases $-1$ and 1: namely, the significant potential is considered to be still a double-well, but with the two local minima coinciding with the pure phases. The most common choice is polynomial of even degree, like the case $\Psi(s)=\frac{1}{4}(s^2-1)^2$.  However, in the case of polynomial potentials, it is worth recalling that it is not
possible to guarantee the existence of physical solutions, that is, solutions for which $-1 \leq  \phi  ( x, t) \leq  1$.
Following, e.g., \cite{Turkinoski}, we get a differential description of the phenomenon of the phase separation as 
\begin{equation}
\partial_t \phi + \text{\text{div}}\ \textbf{J}=0\ \ \text{in }\Omega \times (0,T),
\label{transport}
\end{equation}
where $\phi$ is the order parameter and \textbf{J} is the diffusional flux given by Fick's law,
\begin{equation*}
\textbf{J}=-M(\phi)\nabla \dfrac{\delta \mathcal{U}(\phi)}{\delta \phi}=-M(\phi)\nabla (-\epsilon \Delta \phi+ \Psi'(\phi)),
\end{equation*}
where $\dfrac{\delta \mathcal{U}(\phi)}{\delta \phi}$ is the variational derivative of $\mathcal{U}(\phi)$.
The function $M(\phi)$ is the mobility of the substances and in this work will be considered as a unitary constant (see, for instance, \cite{CM} and \cite{13g} for an analysis of the case of non constant and degenerate mobility, i.e., vanishing at the pure phases). The Cahn-Hilliard equation with constant mobility then reads:
\begin{equation}
\begin{cases}
\partial_t \phi=\Delta\mu\ \ \ \text{in } \Omega\times (0,T),\\
\mu=-\epsilon \Delta \phi+ \Psi'(\phi)\ \ \ \ \ \text{in } \Omega\times (0,T),
\end{cases}
\label{CHH}
\end{equation}
with the initial condition $\phi_0$ and two boundary conditions which are generally the following:
\begin{equation}
\partial_\textbf{n} \phi=0,\qquad \partial_\textbf{n}\mu=0,\qquad  \text{ on } \partial \Omega\times (0,T),
\label{BC}
\end{equation}
with \textbf{n} as the outer normal vector. The former
condition means that no mass flux occurs at the boundary, while the latter requires the
interface to be orthogonal at the boundary.

It is worth noticing that the free energy $\mathcal{U}$ in \eqref{ener_loc} only focuses on short range interactions between particles. Indeed, the gradient square term accounts for the fact that the local interaction energy is spatially dependent and varies across the interfacial surface due to spatial inhomogeneities in the concentration. 
Going back to the general approach of statistical mechanics, the mutual short and long range interactions between particles is described through convolution integrals weighted by interactions kernels. Following this approach, Giacomin and Lebowitz (\cite{giacomin, giacomin1, giacomin2}) observed that a physically more rigorous derivation leads to nonlocal dynamics, which is the nonlocal Cahn-Hilliard equation. In particular, this equation is rigorously justified as a macroscopic limit of microscopic phase segregation models with particles conserving dynamics. In this case, the gradient term is replaced by a nonlocal spatial interaction integral, namely, the energy is defined as
\begin{align}
{\mathcal{E}}(\phi):=-\frac{1}{2}\int_\Omega\int_\Omega J(x-y)\phi(x)\phi(y)\, \d x  \d y+\int_\Omega F(\phi(x)) \, \d x,
\label{en}
\end{align}
where $J$ is a sufficiently smooth symmetric interaction kernel. Note that this functional is characterized by a competition between the mixing entropy $F$ and a nonlocal demixing term. As shown in \cite{giacomin1} (see also \cite{GGG0, GGG} and the references therein), the energy $\mathcal{U}$ can be seen as an
approximation of $\mathcal{E}$, as long as we suitably redefine $F$ as $\widetilde{F}(x,s)=F(s)-\frac{1}{2}(J\ast1)(x)s^2$. In particular, we can rewrite $\mathcal{E}$ as 
\begin{align*}
&\mathcal{E}(\phi)=\frac{1}{4}\int_\Omega\int_\Omega J(x-y)\vert \phi(y)-\phi(x)\vert^2\, \d x  \d y+\int_\Omega \left(F(\phi(x))-\frac{a(x)}{2}\phi^2(x)\right)dx\\&=\frac{1}{4}\int_\Omega\int_\Omega J(x-y)\vert \phi(y)-\phi(x)\vert^2\, \d x  \d y+\int_\Omega \widetilde{F}(\phi(x))dx,
\end{align*}
with $a(x)=(J\ast 1)(x)$. If we formally interpret $\widetilde{F}$ as the potential $\Psi$ of \eqref{ener_loc}, we realize that the (formal) first approximation of the nonlocal interaction is $\frac{k}{2}\vert\nabla\phi\vert^2$, for some $k>0$, as long as $J$ is sufficiently peaked around $0$. In the case $\Omega=\mathbb{T}^3$ (see, e.g., \cite{giacomin2}), the term $J\ast 1$ is a constant: thus $\E$ and $\mathcal{U}$ appear to be very similar. In particular, in this case, corresponding to set $a(x)=\alpha_0$, nonlocal-to-local asymptotics results have been obtained in \cite{Scarpa1,Scarpa2} (see also \cite{GS}) for the nonlocal equation \eqref{nonloc3d} below: namely, the solution to the nonlocal equation converges, under suitable conditions on the data of the problem, to the weak solution of \eqref{CHH}-\eqref{BC}. 

The resulting nonlocal Cahn-Hilliard equation then reads (see \cite{GGG0,GGG})
 \begin{equation}
 \label{nonloc3d}
 \begin{cases}
\partial_t\phi-\Delta\mu=0\quad\text{in }\Omega\times(0,T),\\
\mu=F^\prime(\phi)-J\ast \phi\quad\text{in }\Omega\times(0,T),\\
\partial_\textbf{n}\mu=0\quad\text{ on }\partial\Omega\times(0,T),\\
\phi(\cdot,0)=\phi_0\quad \text{ in }\Omega.
 \end{cases}
 \end{equation}
 From now on we will refer to problem \eqref{CHH}-\eqref{BC} as the local Cahn-Hilliard equation, in order to distinguish it from the nonlocal one in \eqref{nonloc3d}.
 \\The well-posedness theory of Cahn-Hilliard equations with logarithmic (or singular) potential has been studied by many authors in the literature. The local Cahn-Hilliard equation \eqref{CHH}-\eqref{BC} has been studied in \cite{Abels1,Debouche,ElliotL,GGM,Lond,MZ} (see also \cite{CM,GGG} for a review and an insight analysis about this topic).
 Concerning the nonlocal Cahn-Hilliard equation, the physical relevance of nonlocal interactions was already pointed out in the pioneering paper
\cite{vanderwaals} (see also \cite[4.2]{emmerich} and references therein) and studied for different kind of evolution equations, mainly Cahn-Hilliard and phase-field systems (see, e.g., \cite{gillette,colli,zach,  GalGrasselli, rocca}). 
In particular, regarding the nonlocal system \eqref{nonloc3d}, the existence of weak solutions and their uniqueness, and the
existence of the connected global attractor were proven in \cite{Frigeri,FG, FG1}. Moreover, well-posedness and regularity of weak solutions are studied in \cite{GGG0}, namely, in this work the authors establish the validity of the strict separation property in dimension two for the nonlocal Cahn-Hilliard equation \eqref{nonloc3d} with constant mobility and singular potential. This means that if the initial state is not a pure
phase (i.e., $\phi_0\equiv1$ or $\phi_0\equiv-1$), then the corresponding solution stays away from the pure states in finite time, uniformly with respect to the initial datum. Exploiting this crucial property in dimension two, the authors derive straightforward consequences, such as further regularity
results as well as the existence of regular finite dimensional attractors and the convergence of a weak solution to a single equilibrium point. In the recent contribution \cite{GGG}, the same authors propose an alternative argument to prove the strict separation property in dimension two, relying on a De Giorgi's iteration scheme (see \cite[Theorem 4.1]{GGG}). 

In the present work we extend the results of \cite{GGG} to the case of three-dimensional bounded domains, namely we prove the validity of the instantaneous strict separation property in dimension three for the system \eqref{nonloc3d} with singular potential $F$. Our main result is the following: given a weak solution to \eqref{nonloc3d},
\begin{equation}
\forall \tau>0 \ \ \exists \ \delta>0: \vert \phi(x,t)\vert\leq 1-\delta,\quad \text{for a.e. } (x,t)\in \Omega\times(\tau,+\infty),
\label{ta}
\end{equation}
where $\delta$ depends on the parameters of the problem, the initial datum $\phi_0$ and $\tau$. Furthermore, we show that, if the initial datum $\phi_0$ is more regular and already strictly separated from the pure phases, then \eqref{ta} also holds with $\tau=0$, i.e., the solution is uniformly strictly separated at almost any time $t\geq0$. To assess the importance of property \eqref{ta}, similarly to \cite{GGG0}, we infer some additional regularization results for any weak solution and we prove that each weak solution converges to a single stationary state.   

As far as we are aware, this is the first time the instantaneous strict separation property is shown in three-dimensional bounded domains for the Cahn-Hilliard equation with constant mobility and singular (logarithmic) potential. Indeed, the only available result in dimension three regards the nonlocal Cahn-Hilliard equation with degenerate mobility and singular potential and it has been shown in \cite{londen}. For the local Cahn-Hilliard equation the instantaneous separation property has been first proven to hold in \cite{MZ}, but only in dimension two. Concerning dimension three, only the asymptotic (i.e., from some positive time on, depending on the specific initial datum) separation property has been proven in \cite{Abels1} for the local Cahn-Hilliard equation, but nothing is known about its instantaneous (i.e., from \textit{any} positive time on) counterpart. The main issue which so far seemed to be hard to overcome in dimension three for both local and nonlocal cases is the use of the Trudinger-Moser inequality (see, e.g., \cite{Trudi}), which, in dimension $d=2,3$, reads 
\begin{align}
\int_\Omega e^{\vert f(x)\vert}dx \leq Ce^{C\Vert f\Vert^d_{W^{1,d}(\Omega)}},\qquad \forall f\in W^{1,d}(\Omega),
\label{trudi}
\end{align} 
for some positive constant $C$ independent of $f$, but depending on the dimension $d$ and on the Lebesgue $d$-dimensional measure of $\Omega$.
In dimension two this inequality is easy to be handled, since it concerns only the $H^1(\Omega)$ norm of $f$. Indeed, if one assumes that   
\begin{align}
\Fsec(s)\leq Ce^{C\vert F^\prime(s)\vert},\qquad \forall s\in(-1,1),
\label{exp}
\end{align}
for some constant $C>0$ (see, e.g., \cite[(E2)]{GGG} or \cite{GGG0}), which is satisfied by the logarithmic potential 
\begin{equation}
F(s)=\frac{\bar{\alpha}}{2}((1+s)\text{ln}(1+s)+(1-s)\text{ln}(1-s)),\ \ \ \forall s\in[-1,1],
\label{potential}
\end{equation}
then, exploiting \eqref{trudi} as done in \cite{GGG0} or adopting an argument as in \cite[Theorem 3.1]{GGG}, one can control the quantity $\Vert \Fsec(\phi(t))\Vert_{L^p(\Omega)}$, for any $p\geq2$, uniformly in time and this is the key tool to prove the validity of the separation property in two dimensions for example of the nonlocal Cahn-Hilliard equation with constant mobility and singular potential. In the case of three-dimensional bounded domains, \eqref{trudi} leads to the necessity of a control of the $W^{1,3}(\Omega)$ norm of $f$ and this does not seem to be feasible in this context. Thus the proof proposed in \cite{GGG0} does not hold in dimension three. Moreover, also the alternative proof in \cite{GGG} to allow the control of $\Vert \Fsec(\phi(t))\Vert_{L^p(\Omega)}$ is not viable in dimension three, due to the fact that the embedding $H^1(\Omega)\hookrightarrow L^q(\Omega)$ holds only for $q\in[2,6]$, so that a result like \cite[(3.3)-(3.6)]{GGG} cannot be obtained.

 Here we are able to establish the (strict) separation property in three dimensions by avoiding the control of the quantity $\Fsec(\phi(t))$ in any $L^p(\Omega)$ space. We do not assume condition \eqref{exp} on $F$  any more (see assumptions \ref{h3}-\ref{h4} and Remark \ref{esp} below), but we only rely on some natural growth conditions of $F^\prime$ and $\Fsec$ near the endpoints $\pm1$. The idea is to perform a De Giorgi's iteration scheme on each interval of the form $(T-\widetilde{\tau},T)$, with $T>0$ arbitrary and $\widetilde{\tau}$ suitably chosen, similarly to the proof of \cite[Theorem 4.1]{GGG}, but modifying the argument in order to fully exploit the property that $\Fsec(1-2\delta)^{-4}=O(\delta^{4})$, for $\delta>0$ sufficiently small (see \eqref{last}). This is possible in the estimates by treating in a suitable way all the terms leading to the presence of a quantity of the kind $\Fsec(1-2\delta)^{-\gamma}$, with $0\leq\gamma<4$ (see, e.g., the term $Z_2$ in the proof of \cite[Theorem 4.1]{GGG}). To this aim, we first show the validity of a novel Poincar\'{e}-type inequality (Lemma \ref{Poin}), which is applied to a particular family of truncated functions obtained from the weak solution $\phi$ (namely, a family $\phi_{\rx{\rho}}=(\phi-\rx{\rho})^+$, for some suitable $\rx{\rho}\in(0,1)$). This can be obtained heavily relying on the conservation of total mass (i.e., $$\int_\Omega\phi_0(x)dx=\int_\Omega\phi(x,t)dx,$$ for any $t\geq0$), that is one of the most important properties of the solution. By means of this Poincar\'{e}-type inequality, in the De Giorgi's scheme we get, at the end of the estimates, a term of the kind $\Fsec(1-2\delta)^{-4}\delta^{-5}=O(\delta^{-1})$ and this, together with the use of the growth condition of $F^\prime$ near $1$, permits to obtain the strict separation property by choosing a suitably small $\widetilde{\tau}$ depending on $\delta$. Since the size of $\delta$ and the related quantity $\widetilde{\tau}$ do not depend on $T$, we repeat the same argument on each time interval $(T-\widetilde{\tau},T)$ for arbitrary $T>0$, extending the result of the separation property on the entire interval $(\tau,+\infty)$, for $\tau>0$ arbitrarily fixed at the beginning, completing in this way the proof of the validity of \eqref{ta}.\\ As future work, it is worth noticing that the strict separation property could pave the way for the study of other related problems with logarithmic potential in dimension three. For example, one could study the nonlocal Cahn-Hilliard-Oono equation (see, e.g., \cite{DP}), the nonlocal Cahn-Hilliard-Hele-Shaw system (see, e.g., \cite{DellaPorta}) as well as other hydrodynamic phase-field models for binary fluid mixtures of incompressible viscous fluids (see also Remark \ref{extension}).

The paper is organized as follows. In Section \ref{prel} we introduce the functional setting. Section \ref{prel2} is devoted to the presentation some preliminaries, which are essential in the proofs, in particular the new Poincar\'{e}-type inequality. In the same Section we also recall some already-known results concerning well-posedness of the nonlocal Cahn-Hilliard equation and we present a Lemma on geometric convergence of numerical sequences, which is a key tool for De Giorgi's type arguments. Section \ref{main_results} contains the main result concerning the strict separation property in dimension three for the system \eqref{nonloc3d}, together with its proof. In conclusion, in Section \ref{consequences} we present some consequences of the validity of the strict separation property, namely we show some regularization results and we prove that any weak solution to \eqref{nonloc3d} converges to a single equilibrium. 
\section{Mathematical Setting}
\label{prel}

Let $\Omega$ be a smooth bounded domain in $\mathbb{R}^3$. 
The Sobolev spaces are denoted as usual by $W^{k,p}(\Omega)$, where $k\in \mathbb{N}$ and $1\leq p \leq \infty$, with norm $\| \cdot \|_{W^{k,p}(\Omega)}$. The Hilbert space $W^{k,2}(\Omega)$ is denoted by $H^k(\Omega)$ with norm $\|\cdot \|_{H^k(\Omega)}$. In particular, we will adopt the notation
\begin{equation*}
H=L^2(\Omega),\quad V=H^1(\Omega),\quad  V_2=\{v\in H^2(\Omega):\ \partial_\textbf{n}v=0\ \text{on }\partial\Omega\}.
\end{equation*}
Moreover, given a space $X$, we denote by $\textbf{X}$ the space of vectors of three components, each one belonging to $X$.
We then denote by $(\cdot, \cdot)$ the inner product in  $H$ and by $\Vert   \cdot\Vert$ the induced norm. We indicate by $(\cdot,\cdot)_V$ and $\Vert   \cdot\Vert_V$ the canonical inner product and its induced norm in $V$, respectively. We also define the integral mean of a function $f$ as  
$$
\overline{f}:=\dfrac{\int_\Omega f(x)dx}{\vert\Omega\vert},
$$
where $\vert\Omega\vert$ stands for the three-dimensional Lebesgue measure of the set $\Omega$.
We then introduce
$$
H_0= \lbrace v \in H: \overline{f}=0 \rbrace, \quad 
V_0= \lbrace v \in V: \overline{f}=0 \rbrace, \quad 
V_0' = \lbrace v \in V': \frac{\l f,1 \r}{|\Omega|}=0 \rbrace,
$$
endowed with the norms of $H$, $V$ and $V'$. Thanks to the Poincar\'{e}-Wirtinger inequality, it follows that $(\| \nabla u\|_{L^2(\Omega)}^2+ |\overline{u}|^2)^\frac12$ is a norm on $V$ equivalent to $\| u\|_V$. 
The Laplace operator $A_0: V_0 \rightarrow V_0'$ defined by $<A_0 u,v>=(\nabla u,\nabla v)$ is an isomorphism. We denote by $\mathcal{N}$ its inverse map and we set $\Vert  f\Vert  _*:=\Vert  \nabla \mathcal{N}f\Vert$, which is a norm on $V_0'$ equivalent to the canonical one. Moreover, we recall that 
\begin{align}
\label{equiv}
\Vert f-\overline{f}\Vert_*^2+\vert\overline{f}\vert^2
\end{align}
is a norm $V^\prime$ which is equivalent to the standard one. 
Next, we recall the following Gagliardo-Nirenberg's inequality (see, e.g., \cite[Ch.9]{Brezis})
\begin{align}
\label{Gagl}
&\| u\|_{L^p(\Omega)}\leq C(p) \|u\|^{\frac{6-p}{2p}}\|u\|_{V}^{\frac{3(p-2)}{2p}}, \quad &&\forall \, u \in V,\qquad \forall p\in[2,6], 
\end{align}
where the constant $C(p)$ depends on $\Omega$ and $p$. From this inequality, in the case $p=\frac{10}{3}$ we get
\begin{align}
\label{sobolev}
&\| u\|_{L^\frac{10}{3}(\Omega)}\leq \widehat{C} \|u\|^{\frac{2}{5}}\|u\|_{V}^{\frac{3}{5}}, \quad &&\forall \, u \in V, 
\end{align}
with $\widehat{C}>0$ depending on $\Omega$. 

\noindent

\section{Preliminaries}
\label{prel2}
Here we present some preliminary results, which are essential for the proof of our main theorem.
\subsection{A Poincar\'{e}-type inequality}
 First we state the following generalized version of the well known Poincar\'{e}'s inequality:
\begin{lemma}
	\label{Poin}
	Let $I$ be either a compact interval or an interval of the kind $[\tau,+\infty)$, with $\tau>0$. \rx{Let $\mathcal{K}\subset \R$ be a set of indices} and $\{f_{\rx{\rho}}\}_{\rx{\rho\in \mathcal{K}}}\subset L^{\extra{\infty}}(I;V)\cap C(I;H)$. Assume also that, for any \rx{$\rho\in \mathcal{K}$} and for any $t\in I$, $f_{\rx{\rho}}(t)\equiv 0$ on the set $E(t):=\{x\in \Omega: g(t,x)\leq 1-2\delta\}\subset\Omega$, with $g\in C(I;L^q(\Omega))$, $q\geq1$, and $\delta\in(0,\frac{1}{2})$. Moreover, for a fixed $\varepsilon>0$ sufficiently small, assume that for any $t\in I$ the set $\{x\in \Omega: g(t,x)\leq 1-2\delta-\varepsilon\}\subset E(t)$ has strictly positive Lebesgue measure. In the case the interval $I$ is $[\tau,+\infty)$, assume additionally that for any sequence $\{t_l\}_l$, such that $t_l\to \infty$ as $l\to \infty$, there exists a (non-relabeled) subsequence $\{t_l\}_l$, a function $g^\star\in L^r(\Omega)$, $r\geq1$, and $\widetilde{\varepsilon}>0$, 
	such that $g(t_l)\to g^\star$ strongly in $L^r(\Omega)$ as $l\to\infty$ and the set $\{x\in \Omega: g^\star(x)\leq 1-2\delta-\widetilde{\varepsilon}\}$ has strictly positive Lebesgue measure.\\
	Then there exists a uniform (in $\rx{\rho}$ and $t$) constant $C_P>0$ 
	such that 
	\begin{equation}
	\Vert f_{\rx{\rho}}(t)\Vert\leq C_P\Vert\nabla f_{\rx{\rho}}(t)\Vert\qquad \forall t\in I,\qquad \forall \rx{\rho\in \mathcal{K}}.
	\label{eq}
	\end{equation}
\end{lemma}
\begin{remark}
    \extra{Being $\{f_{\rx{\rho}}\}_{\rx{\rho}}\subset C(I;H)\cap L^\infty(I;V)\hookrightarrow C_{w}(I;V)$, where $C_w(I;V)$ denotes the $V$-valued weakly continuous functions (see, e.g., \cite[Lemma II.5.9]{Boyerlibro}), it  makes sense to ask for conditions at \textit{any} time $t\in I$.}
\end{remark}

	\begin{proof}
		Being $\{f_{\rx{\rho}}\}_{\rx{\rho}}\subset {\extra{C_w}}(I;V)$, $f_{\rx{\rho}}(t)\in V$ for any \rx{$\rho\in \mathcal{K}$} and any $t\in I$. Assume by contradiction that \eqref{eq} is false. Then there exist a sequence $\{\rx{\rho}_n\}_{n\in\N}\subset \rx{\mathcal{K}}$ and a sequence $\{t_n\}_{n\in\N}\subset I$ such that
		$$
		\Vert f_{\rx{\rho}_n}(t_n)\Vert> n\Vert\nabla f_{\rx{\rho}_n}(t_n)\Vert,\quad \forall n\in \N.
		$$ 
		We then set $w_n:=\dfrac{f_{\rx{\rho}_n}(t_n)}{\Vert f_{\rx{\rho}_n}(t_n) \Vert}$, with $\Vert w_n\Vert=1$. 
		We need to consider two cases:
		\begin{enumerate}
		    \item Either the interval $I$ is compact or there exists a non-relabeled subsequence of $\{t_{n}\}_n$ which is entirely contained in the set $[\tau, M]\subset I$, for some  $M<+\infty$.
		    In this case there exists another non-relabeled subsequence of times and $t^\star\in I$, with $t^\star<+\infty$, such that $t_n\to t^\star$.
		    
		 Now notice that, being $g\in C(I;L^q(\Omega))$, $q\geq1$, we get $g(t_n)\to g(t^\star)$ in $L^q(\Omega)$. Therefore, there exists a subsequence $\{g(t_{n_j})\}_j$ such that, as $j\to\infty$, 
		 $$
		 g(t_{n_j})\to g(t^\star)\text{ a.e. in }\Omega.
		 $$
		 Let us now set $D:=\{x\in \Omega: g(t^\star,x)\leq 1-2\delta-\varepsilon\}$, and $$\alpha=\vert D\vert>0,$$
		 which is possible by assumption.
		 Then by Severini-Egorov Theorem (notice that $\Omega$ has finite measure, so this theorem can be applied), there exists a measurable subset $B\subset\Omega$ such that $\vert B\vert <\frac{\alpha}{2}$ and such that, as $j\to\infty$,
		 $$
		 g(t_{n_j})\to g(t^\star)\text{ uniformly on } \Omega\setminus B.
		 $$ 
		 Therefore, we also deduce that $\vert D\setminus B\vert>\frac{\alpha}{2}>0$ and that also 
		 
		 $$
		 g(t_{n_j})\to g(t^\star)\text{ uniformly on }   D\setminus B.
		 $$
		 This means that there exists a $\overline{J}\in \N$ such that, for any $x\in D\setminus B$,
		 $$
		 \vert g(t_{n_j},x)-g(t^\star, x)\vert<\varepsilon\qquad \forall j\geq \overline{J},
		 $$
		 implying that, for any $x\in D\setminus B$, by definition of the set $D$,
		 $$
		 g(t_{n_j},x)=g(t_{n_j},x)-g(t^\star, x)+g(t^\star, x)\leq \varepsilon+1-2\delta-\varepsilon=1-2\delta,\quad  \forall j\geq \overline{J} 
		 $$
		 This means, by the assumptions, that 
		 $$
		 D\setminus B \subset E(t_{n_j})\subset \{x\in\Omega:\ w_{n_j}(x)=0\} \quad \forall j\geq \overline{J},
		 $$
		 implying 
		 $$
		 D\setminus B\subset \bigcap_{j\geq \overline{J}}\{x\in\Omega:\ w_{n_j}(x)=0\},\qquad \vert D\setminus B\vert>\frac{\alpha}{2}.
		 $$
		 \label{a}
		 \item \label{b}
		  The interval $I$ is of the form $[\tau,+\infty)$ and there are no bounded subsequences of $\{t_n\}_n$, i.e. $t_n\to+\infty$ as $n\to \infty$. In this case we have by assumption that, up to a non-relabeled subsequence, there exists $g^\star\in L^r(\Omega)$, $r\geq1$, such that $g(t_n)\to g^\star$ strongly in $L^r(\Omega)$. Thus there exists a subsequence $\{g(t_{n_j})\}_j$ such that 
		 $$
		 g(t_{n_j})\to g^\star\text{ a.e. in }\Omega.
		 $$
		 As in case \eqref{a}, we set $D:=\{x\in \Omega: g^\star(x)\leq 1-2\delta-\widetilde{\varepsilon}\}$, and $$\alpha=\vert D\vert>0,$$
		 which is again possible by assumption. Then we can repeat exactly the same arguments as in case \eqref{a} to obtain again that 
			$$
		 D\setminus B \subset E(t_{n_j})\subset \{x\in\Omega:\ w_{n_j}(x)=0\} \quad \forall j\geq \overline{J},
		 $$
		 implying 
		 $$
		 D\setminus B\subset \bigcap_{j\geq \overline{J}}\{x\in\Omega:\ w_{n_j}(x)=0\},\qquad \vert D\setminus B\vert>\frac{\alpha}{2}.
		 $$
		 Clearly notice that in this case the set $B$ will be such that there exists a $\overline{J}\in \N$ such that,  for any $x\in D\setminus B$,
		 $$
		 \vert g(t_{n_j},x)-g^\star( x)\vert<\widetilde{\varepsilon}\qquad \forall j\geq \overline{J}.
		 $$
		\end{enumerate}
In both cases \eqref{a} and \eqref{b}, being $w_{n_j}$ uniformly bounded in $V$, there exists $w\in V$ such that, by Rellich-Kondrachov Theorem, as $j\to\infty$,
		$$
		w_{n_j}\rightharpoonup w\text{  in } V,\quad w_{n_j}\to w\text{ in }H,\quad \nabla w_{n_j}\rightharpoonup \nabla w\text{ in }H,
		$$
		up to a non-relabeled subsequence.
		Moreover, being $\Vert\nabla w_{n_j}\Vert< \frac{1}{{n_j}}$, we deduce, by weak lower sequential semicontinuity of the $L^2$-norm, that $\nabla w\equiv 0$ almost everywhere in $\Omega$ and thus, being $\Omega$ connected, $w\equiv \kappa$ almost everywhere in $\Omega$, with $\kappa$ constant.
		 Therefore, since also, up to another subsequence, $w_{n_j}\to w$ almost everywhere in $\Omega$, we have $w\equiv 0$ on $D\setminus B$ (of positive Lebesgue measure) up to a zero measure set. But this clearly implies that $\kappa=0$, which is a contradiction, since $\Vert w\Vert=1$, being $\Vert w_{n_j}\Vert=1$ and $w_{n_j}\to w$ in $H$ as $j\to\infty$. This concludes the proof.
	\end{proof}	
\subsection{The state of the art for the three-dimensional nonlocal Cahn-Hilliard equation}
For the sake of completeness we state here the already-known results concerning the nonlocal Cahn-Hilliard equation with constant mobility and singular potential in three dimensional bounded domains. We first consider the following assumptions: 

\begin{enumerate}[label=(\subscript{H}{{\arabic*}})]
	\item $J \in W_{loc}^{1,1}(\R^3)$, with $J({x})=J(-{x})$.
	\label{h1}
	\item $F\in C([-1, 1]) \cap C^2(-1,1)$ fulfills
	\begin{equation*}
	\lim_{s\to-1} F'(s)=-\infty, \quad \lim_{s\to1} F'(s)=+\infty,\quad F''(s)\geq{\alpha\extra{>0}}, \quad\forall\ s\in(-1,1).
	\end{equation*}
	We extend $F(s)=+\infty$  for any $s\notin[-1, 1]$.  Without loss of generality, $F(0) = 0$ and $F'(0)=0$. In particular, this entails that $F (s) \geq 0$ for any $s \in[-1, 1]$. Also, we assume that there exists $\gamma\in(0,1)$ such that $F''$ is nondecreasing in $[1-\gamma,1)$ and nonincreasing in $(-1,-1+\gamma]$.
	\label{h3}
\end{enumerate}
 We then have the following 
\begin{theorem}
	\label{known}
	Assume that \ref{h1}-\ref{h3} hold and also that $\phi_0\in L^\infty(\Omega)$ such that $\Vert\phi_0\Vert_{L^\infty}\leq 1$ and $\vert\overline{\phi}_0\vert=m<1$. Then there exists a unique weak solution to \eqref{nonloc3d} such that, for any $T>0$,
	\begin{align*}
	&\phi\in L^\infty(\Omega\times(0,T)):\quad \forall t> 0, \quad \vert\phi(t)\vert<1,\quad\text{a.e. in }\Omega,\\&
	\phi\in L^2(0,T;V)\cap H^{1}(0,T;H),\\&
	\mu\in L^2(0,T;V),\qquad F^\prime(\phi)\in L^2(0,T;V),
	\end{align*}
	such that 
	\begin{align}
	&<\partial_t\phi,v>+(\nabla\mu,\nabla v)=0\quad\forall v \in V,\qquad \text{a.e. in }(0,T),\label{phi}\\&
	\mu=F^\prime(\phi)-J\ast\phi\quad\text{a.e. in }\Omega\times(0,T),
	\label{mu}
	\end{align}
	and $\phi(\cdot,0)=\phi_0(\cdot)$ in $\Omega$.
	The weak solution also satisfies the energy identity ($\mathcal{E}$ is defined in \eqref{en})
	\begin{align}
	\mathcal{E}(\phi(t))+\int_s^{t}\Vert\nabla\mu(\tau)\Vert^2d\tau= \mathcal{E}(\phi(s)),\qquad \forall\ 0\leq s\leq t<\infty.
	\label{dissipative}
	\end{align}
	 Moreover, for any $\tau>0$, 
	\begin{align}
	&\label{dt}\sup_{t\geq \tau}\Vert\partial_t\phi(t)\Vert_{V^\prime}+\sup_{t\geq \tau}\Vert\partial_t\phi\Vert_{L^2(t,t+1,H)}\leq \frac{K_0}{\sqrt{\tau}},\\&\label{H1}
	\sup_{t\geq \tau}\Vert\mu(t)\Vert_{V}+\sup_{t\geq \tau}\Vert\phi(t)\Vert_{V}\leq \frac{K_0}{\sqrt{\tau}},\\&
	\Vert F^\prime(\phi)\Vert_{L^\infty(\tau,t;V)}+\Vert\mu\Vert_{L^2(t,t+1,V_2)}\leq K_1,\quad \forall t\geq \tau,
	\label{Fprime}	
	\\&\Vert\nabla\mu\Vert_{L^q(t,t+1;L^p(\Omega))}+\Vert\nabla\phi\Vert_{L^q(t,t+1;L^p(\Omega))}\leq K_2,\qquad \text{if }\quad \frac{3p-6}{2p}=\frac{2}{q},\quad \forall p\in[2,6],\quad \forall t\geq \tau,	
	\label{Lp}
	\end{align}
	where the positive constant $K_0$ depends only on the initial datum energy $\E(\phi_0)$, $\overline{\phi}_0$, $\Omega$ and the parameters of the system, whereas $K_1=K_1(\tau)$ and $K_2=K_2(\tau)$ also depend on $\tau$. Furthermore $K_2$ depends on also $q,p$.
	In conclusion, it holds the following continuous dependence estimate: for every two weak solutions $\phi_1$ and $\phi_2$ to \eqref{nonloc3d} on $[0,T]$, with initial data $\phi_{01}$ and $\phi_{02}$, respectively, we have, for all $t\in[0,T]$,
	$$
	\Vert\phi_1(t)-\phi_2(t)\Vert_{V^\prime}^2\leq \Vert\phi_{01}-\phi_{02}\Vert_{V^\prime}^2+K\vert \overline{\phi}_{01}-\overline{\phi}_{02}\vert e^{CT},
	$$
	where $C$ is a positive constant and
	$$
	K=C\left(\Vert F^\prime(\phi_1)\Vert_{L^1(0,T;{L^1(\Omega)})}+\Vert F^\prime(\phi_2)\Vert_{L^1(0,T;{L^1(\Omega)})}\right).
	$$
\end{theorem}
\rx{\begin{remark}
The proof of the above theorem can be found in \cite[Theorems 3.4, 4.1, Proposition 4.2]{GGG0} and \cite[Proposition 3.1]{DellaPorta}, see also \cite[Theorem 4.1]{AGGP} and \cite[Theorem 2.2]{PS} for a comprehensive result in the more general case of an advective nonlocal Cahn-Hilliard equation in two and three dimensions, respectively. In particular, we refer to \cite[Theorem 4.1, (4.4)]{AGGP} and \cite[Proposition 3.1, (3.53)]{DellaPorta}, which still hold in the non-advective case $\uu=\textbf{0}$, for the validity of the energy identity \eqref{dissipative}, whereas \eqref{dt} is shown in \cite[Theorem 4.1, (4.2)]{GGG0}. Estimates \eqref{H1}-\eqref{Fprime} can be found in \cite[Theorem 4.1, (4.3)]{GGG0}-\cite[Proposition 4.2, (4.7)]{GGG0}, while \eqref{Lp} is shown in \cite[Proposition 4.2, (4.9)]{GGG0}.
\end{remark}}
\begin{remark}
	If we assume additionally that $\nabla F^\prime(\phi_0)\in \textbf{H}$ we can actually extend \eqref{dt}-\eqref{Lp} to $\tau=0$, since the initial datum is more regular and one can argue as in \cite[Sec.4]{DellaPorta} to obtain the desired regularity departing from the initial time. This means that the solution $\phi$ with initial datum $\phi_0$ is indeed a strong solution to problem \eqref{nonloc3d}.
		\label{strong}
\end{remark}
\begin{remark}
	Notice that from condition \eqref{Fprime} we can also deduce by Sobolev embeddings that 
		\begin{align}
	\Vert F^\prime(\phi)\Vert_{L^\infty(\tau,\infty;L^p(\Omega))}\leq K_3(\tau,p),\qquad \forall p\in[1,6],
	\label{Fprime2}		
	\end{align}
	where $K_3(\tau,p)$ depends on $K_1$, $\Omega$ and $p$.
\end{remark}
\begin{remark}
	We highlight that the previous theorem and our following main result concerning the strict separation property in dimension three heavily rely on the assumption $\overline{\phi}_0\in (-1, 1)$ (see also \cite{Kenmochi} for the local Cahn-Hilliard equation). This is physically reasonable
	since $\overline{\phi}_0=1$ (or $\overline{\phi}_0=-1$) means that the initial condition is a pure phase, so that no
	phase separation takes place in $\Omega$, unless we assume the existence of a source or reaction term (see, for instance \cite{GSS}).
\end{remark}
\subsection{A lemma on geometric convergence of sequences}
We present here one of the key tools for the application of De Giorgi's iteration argument. This Lemma can be found, e.g., in \cite[Ch. I, Lemma 4.1]{DiBenedetto}, \cite[Ch.2, Lemma 5.6]{Lady}, and it has also been proposed in \cite[Lemma 4.3]{GGG}.
\begin{lemma}
	\label{conv}
	Let $\{y_n\}_{n\in\N\cup \{0\}}\subset \R^+$ satisfy the recursive inequalities 
	\begin{align}
	y_{n+1}\leq Cb^ny_n^{1+\varepsilon},
	\label{ineq}\qquad \forall n\geq 0,
	\end{align}
	for some $C>0$, $b>1$ and $\varepsilon>0$. If 
	\begin{align}
	\label{condition}
	y_0\leq \theta:= C^{-\frac{1}{\varepsilon}}b^{-\frac{1}{\varepsilon^2}}, 
	\end{align}
	then 
	\begin{align}
	y_n\leq \theta b^{-\frac{n}{\varepsilon}},\qquad \forall n\geq 0,
	\label{yn}
	\end{align}
	and consequently $y_n\to 0$ for $n\to \infty$.
\end{lemma}
\begin{proof}
	The proof can be easily carried out directly by induction. Indeed, the case $n=0$ is trivial. Then assume that \eqref{yn} holds for $n$. We prove that it also holds for $n+1$. In particular we have by \eqref{ineq} and recalling \eqref{condition},
	$$
	y_{n+1}\leq Cb^{n}y_n^{1+\varepsilon}\leq Cb^n\theta^{1+\varepsilon}b^{{-\frac n\varepsilon}(1+\varepsilon)}=C\theta^{1+\varepsilon}b^{-\frac n\varepsilon}=\theta b^{-\frac{n+1}{\varepsilon}}C\theta^\varepsilon b^{\frac 1\varepsilon}\leq \theta b^{-\frac{n+1}{\varepsilon}},
	$$
	where we exploited the definition of $\theta$ in \eqref{condition}. This means that \eqref{yn} also holds for $n+1$, concluding the proof by induction.
\end{proof}

We now present our main results, concerning the instantaneous strict separation property in three dimensional bounded domains.
\section{Main results}	
\label{main_results}	
Let us assume, additionally to \ref{h3}, the following hypotheses on the singular potential $F$:
\begin{enumerate}[label=(\subscript{H}{{3}})]
	\item As $\delta\to 0^+$ we assume
	\begin{align}
	\frac{1}{F^{\prime}(1-2\delta)}=O\left(\frac{1}{\vert\ln(\delta)\vert}\right),\quad \frac{1}{F^{\prime\prime}(1-2\delta)}=O(\delta),
	\label{F}
	\end{align}
	and analogously
	\begin{align}
	\dfrac{1}{\vert F^{\prime}(-1+2\delta)\vert }=O\left(\frac{1}{\vert\ln(\delta)\vert}\right),\qquad \dfrac{1}{\Fsec(-1+2\delta)}=O\left(\delta\right).
	\label{F2}
	\end{align}
	\label{h4}
\end{enumerate}
 \begin{remark}
 Notice that these conditions are verified by the logarithmic potential \eqref{potential}. Indeed, it holds 
 \begin{align*}
 F^\prime(s)=\frac{\overline{\alpha}}{2}\ln\left(\dfrac{1+s}{1-s}\right),\qquad \Fsec(s)=\dfrac{\overline{\alpha}}{1-s^2},
 \end{align*} 
 thus 
 \begin{align*}
 F^\prime(1-2\delta)=\frac{\overline{\alpha}}{2}\ln\left(\dfrac{1-\delta}{\delta}\right),\qquad \Fsec(1-2\delta)=\dfrac{\overline{\alpha}}{4\delta(1-\delta)},
 \end{align*} 
 and 
  \begin{align*}
 F^\prime(-1+2\delta)=\frac{\overline{\alpha}}{2}\ln\left(\dfrac{\delta}{1-\delta}\right),\qquad \Fsec(-1+2\delta)=\dfrac{\overline{\alpha}}{4\delta(1-\delta)},
 \end{align*} 
 clearly implying assumption \ref{h4}.
 \end{remark}
\begin{remark}
As already pointed out in the Introduction, assumption \ref{h4} does not make any explicit reference to the typical extra condition \eqref{exp}.  Indeed, as far as we know, this is the first proof of the instantaneous separation property concerning nonlocal Cahn-Hilliard equation with constant mobility and singular potential (problem \eqref{nonloc3d}) in which it is not exploited any constraint on $\Vert \Fsec(\phi(t))\Vert_{L^q(\Omega)}$, for some $q\geq2$ and for almost any $t\geq \tau$, with $\tau>0$. 
\label{esp}
Indeed, in our proof we simply rely on some natural growth conditions of $F^\prime$ and $\Fsec$ near the endpoints $\pm1$.
Note that assumptions \ref{h3}-\ref{h4} on the potential $F$ are somehow minimal, in the sense that the proof of the separation property in dimension three works only in this case (or for more singular potentials than the logarithmic one). This seems to suggest that the use of the logarithmic potential when modeling phase separation phenomena with the help of the nonlocal Cahn-Hilliard equation with constant mobility could be a good choice, since it preserves all the basic physical properties expected from the solution.
\end{remark}

We can now state our main
\begin{theorem}
Let $\Omega\subset \R^3$ be a smooth bounded domain and let assumptions \ref{h1}-\ref{h4} hold. Assume that $\phi_0\in L^\infty(\Omega)$ such that $\Vert\phi_0\Vert_{L^\infty}\leq 1$ and $\vert\overline{\phi}_0\vert=m<1$. Then for any $\tau>0$ there exists $\delta\in (0,1)$, depending on $\tau$, $m$ and the initial datum, such that the unique weak solution to problem \eqref{nonloc3d} given in Theorem \ref{known} satisfies
$$
\vert \phi(x,t)\vert\leq 1-\delta,\quad \text{for a.e. } (x,t)\in \Omega\times(\tau,+\infty),
$$
i.e., the instantaneous strict separation property from the pure phases $\pm 1$ holds.
\label{sep}	
\end{theorem}
\begin{remark}
	Observe that the quantity $\delta$ given in the theorem strongly depends on the specific entire trajectory, therefore, by the uniqueness of the solution, on the initial datum $\phi_0$. This means that we cannot have an explicit dependence of $\delta$, e.g., on the initial datum energy.
\end{remark}
As a byproduct of the main theorem, we also prove that, if the initial datum $\phi_0$ is more regular and already separated from the pure phases, i.e., there exists $\delta_0\in(0,1]$ such that 
$$
\Vert\phi_0\Vert_{L^\infty(\Omega)}\leq 1-\delta_0,
$$
then the unique solution $\phi$ departing from $\phi_0$, which is now strong from the time $t=0$ (see Remark \ref{strong}), is strictly separated on $[0,+\infty)$, i.e., it remains separated from the pure phases \textit{uniformly} for almost any time $t\geq 0$. In particular, we have 
\begin{corollary}
	Under the same hypotheses of Theorem \ref{sep}, if we assume additionally that $\nabla F^\prime(\phi_0)\in \textbf{H}$, and that $\phi_0$ is strictly separated, i.e.,  
	there exists $\delta_0\in(0,1]$ such that 
	$$
	\Vert\phi_0\Vert_{L^\infty(\Omega)}\leq 1-\delta_0,
	$$
	 then there exists $\delta\in (0,1)$, depending on $\tau$, $m$, $\delta_0$ and the initial datum, such that the unique  strong solution to problem \eqref{nonloc3d} given in Remark \ref{strong} satisfies
	$$
	\vert \phi(x,t)\vert\leq 1-\delta,\quad \text{for a.e. } (x,t)\in \Omega\times[0,+\infty),
	$$
	i.e., the instantaneous strict separation property from the pure phases $\pm 1$ holds for almost any time $t\geq0$.
	\label{sep2}	
\end{corollary}
\begin{remark}
	\label{proof}
Observe that, since by Theorem \ref{sep} the solution $\phi$ in \rx{Corollary} \ref{sep2} is strictly separated on time sets of the kind $(\tau,+\infty)$, for any $\tau>0$, it is enough to show that there exists an interval $[0,T_1]$ ($T_1>0$) on which the solution is separated to obtain the strict separation over $[0,+\infty)$, choosing $\tau=T_1$. As it will be clear from the proof of \rx{Corollary} \ref{sep2}, $T_1$ can be explicitly computed as a function of the parameters of the problem and the initial datum.
\end{remark}
\subsection{Proof of Theorem \ref{sep}}
We divide the proof into two steps. In the first one we show that we can apply Lemma \ref{Poin} to a specific family of functions, which will be of essential importance in the second step, when we adopt a De Giorgi's iteration scheme (as in \cite[\rx{Theorem} 4.1]{GGG}) to obtain the desired result.
\subsection*{Step 1. Application of Lemma \ref{Poin} to a family of truncated functions}
Let us consider the unique solution $\phi$ departing from $\phi_0$, whose existence and regularity is stated in Theorem \ref{known}.
We make the following observations: first fix any $\tau>0$.
\begin{itemize}
	\item Being $\vert \overline{\phi}_0\vert \leq m<1$, there exists $\widehat{\delta}>0$ and an $\varepsilon
	>0$ such that 
	\begin{align}
	m\leq1-2\widehat{\delta}-\varepsilon.
	\label{m}
	\end{align}
	In particular we may choose $\varepsilon:=\frac{1-m}{2}>0$ and $\widehat{\delta}:=\frac{1-m}{4}>0$.
Thanks to the conservation of total mass, we have that for any \rx{$\rho\in\R^+$, $\rho\geq 1-2\widehat{\delta}$,} and for any $t\in[0,+\infty)$, the function 
	\begin{align}
	\phi_{\rx{\rho}}(x,t):=(\phi(x,t)-\rx{\rho})^+
	\label{phik}
	\end{align}
	vanishes on the set \extra{(independent of $\rho$)} 
	\begin{align}
	E(t):=\{x\in \Omega: \phi(x,t)\leq 1-2\widehat{\delta}\},
	\label{e1}
	\end{align}
	which is such that 
	\begin{align}
	\vert\{x\in \Omega: \phi(x,t)\leq 1-2\widehat{\delta}-\varepsilon\}\vert>0,\qquad \forall t\geq 0.
	\label{e2}
	\end{align}
	\begin{proof}
		To prove this observation, let us assume by contradiction that, for some $\tilde{t}\geq 0$, 
		$$
		\vert\{x\in \Omega: \phi(x,\tilde{t})\leq 1-2\widehat{\delta}-\varepsilon\}\vert=0.
		$$
		By the conservation of total mass we get, for any $t\geq0$, 
		$$
		(1-2\widehat{\delta}-\varepsilon)\vert \Omega\vert\geq m\vert \Omega\vert\geq\int_\Omega \phi_0(x)dx =\int_\Omega \phi(x,t)dx,
		$$
		but then \rx{we get a contradiction, since $\vert \Omega\vert=\vert\{x\in \Omega: \phi(x,\tilde{t})> 1-2\widehat{\delta}-\varepsilon\}\vert$ and
		\begin{align*}
		&(1-2\widehat{\delta}-\varepsilon)\vert \Omega\vert\geq\int_\Omega \phi(x,\tilde{t})dx>(1-2\widehat{\delta}-\varepsilon)\vert\{x\in \Omega: \phi(x,\tilde{t})> 1-2\widehat{\delta}-\varepsilon\}\vert.
		\end{align*}
		}
	\end{proof}
	\item We aim to apply Lemma \ref{Poin} with \rx{$\mathcal{K}=[1-2\widehat{\delta},1]$}, ${\{f_{\rho}\}_{\rx{\rho\in\mathcal{K}}}=\{\phi_{\rx{\rho}}\}_{\rx{\rho\in\mathcal{K}}}}$, $I=[\tau,+\infty)$, $g=\phi$, $\delta=\widehat{\delta}$, $\widetilde{\varepsilon}=\varepsilon$. Indeed we verify all the assumptions:
	\begin{enumerate}[label=(\alph*)]
	    \item We have $\{\phi_{\rx{\rho}}\}_{\rx{\rho}}\subset L^\infty(I;V)\cap C(I;H)$, $\phi\in C(I;H)$, and \eqref{e1} and \eqref{e2} hold for any $t\in I$.
	    \item Let $\{t_l\}_l$ be any sequence such that $t_l\to \infty$. By \eqref{H1}, there exists a constant $C(\tau)>0$ such that 
	    $$
	    \sup_{t\geq \tau}\Vert\phi\Vert_{V}\leq C(\tau).
	    $$
	    Therefore, being $V$ reflexive, there exists a (non-relabeled) subsequence $\{t_l\}_l$ and a function $g^\star\in V$ (which could depend on the subsequence) such that, as $l\to \infty$,
	    $$
	    \phi(t_l)\rightharpoonup g^\star\qquad \text{in } V,
	    $$
	    implying by compactness that
	    \begin{align}
	     \phi(t_l)\to g^\star\qquad \text{in } H.
	    \label{H}
	    \end{align}
	    Now notice that this strong convergence also implies, by the conservation of total mass, that 
	    $$
	    \int_\Omega \phi_0(x)dx =\int_\Omega \phi(x,t_l)dx\to \int_\Omega g^\star(x)dx,
	    $$
	    and thus also $g^\star$ enjoys the same total mass as the initial datum $\phi_0$:
	    $$
	    \int_\Omega g^\star(x)dx=\int_\Omega \phi_0(x)dx.
	    $$
	    This means that we can repeat exactly the same argument as the one adopted to get \eqref{e2} to infer 
	    	\begin{align}
	\vert\{x\in \Omega: g^\star(x)\leq 1-2\widehat{\delta}-\varepsilon\}\vert>0,
	\label{e2bis}
	\end{align}
	so that, having chosen $\widetilde{\varepsilon}=\varepsilon$ and $g=\phi$, thanks to \eqref{H}-\eqref{e2bis}, we have completed the verification of the assumptions of Lemma \ref{Poin}.
	\end{enumerate}
	In the end we can conclude that there exists a uniform (in $\rx{\rho}$ and $t$) constant $C_{P,+}>0$ such that 
	\begin{equation}
	\Vert \phi_{\rx{\rho}}(t)\Vert\leq C_{P,+}\Vert\nabla \phi_{\rx{\rho}}(t)\Vert,
	\label{eq2a}
	\end{equation}
	for any $t\in [\tau,+\infty)$ and any \rx{$\rho\in [1-2\widehat{\delta},1]$}.
	\item Since in the last part of the proof we need to reproduce all the arguments on the functions 
	\begin{align}
	\widetilde{\phi}_{\rx{\rho}}(x,t):=(\phi(x,t)+\rx{\rho})^-=(-\phi(x,t)
	-\rx{\rho})^+,
	\label{phik1}
	\end{align}
	with $\rx{\rho}\geq 1-2\widehat{\delta}$,
	we observe that \eqref{e1} and \eqref{e2} still hold substituting $\phi$ with $-\phi$, simply because, again by the conservation of mass, $m\vert\Omega\vert\geq \int_\Omega -\phi(x,t)dx$ for any $t\geq \tau$. Therefore again the assumptions of Lemma \ref{Poin} are satisfied (with $g=-\phi$), and thus that there exists a uniform (in $\rx{\rho}$ and $t$) constant $C_{P,-}>0$ (which is possibly different from $C_{P,+}$) such that 
	\begin{equation}
	\Vert \widetilde{\phi}_{\rx{\rho}}(t)\Vert\leq C_{P,-}\Vert\nabla \widetilde{\phi}_{\rx{\rho}}(t)\Vert
	\label{eq2bis}
	\end{equation}
	for any $t\in [\tau,+\infty)$ and for any \rx{$\rho\in [1-2\widehat{\delta},1]$}. 
	Therefore, we introduce the constant $C_P:=\max\{C_{P,+},C_{P,-}\}$ so that both \eqref{eq2a} and \eqref{eq2bis} hold with the same constant $C_P$, i.e.,
		\begin{equation}
	\Vert \phi_{\rx{\rho}}(t)\Vert\leq C_{P}\Vert\nabla \phi_{\rx{\rho}}(t)\Vert, \qquad 	\Vert \widetilde{\phi}_{\rx{\rho}}(t)\Vert\leq C_{P}\Vert\nabla \widetilde{\phi}_{\rx{\rho}}(t)\Vert,
	\label{eq2}
	\end{equation}
	 for any $t\geq \tau$ and any \rx{$\rho\in [1-2\widehat{\delta},1]$}.
	 Note that the constant $C_P$ depends on the specific solution $\phi$ we used, thus, since $\phi$ is uniquely determined by $\phi_0$, we have that $C_P$ depends in a nontrivial way on the initial datum.
\end{itemize}
\subsection*{Step 2. De Giorgi's iteration scheme}
We perform a De Giorgi's iteration scheme following the one presented in \cite[Lemma 4.1]{GGG}. Let us fix $\delta$ sufficiently small such that $\delta\leq \widehat{\delta}$, so that \eqref{eq2} holds for any $\rx{\rho\in[1-2\delta,1]}$. Set then $\widetilde{\tau}>0$ such that it holds 
\rx{\begin{align}
\widetilde{\tau}= \dfrac{2^{-20}\delta^5\left(\Fsec(1-2\delta)\right)^{4}F^\prime(1-2\delta)}{3C(\tau)\Vert\nabla J\Vert_{{L^1(B_r)}}^{5}\widehat{C}^{\frac{3}{2}}\left(1+C_P^2\right)^{\frac{3}{2}}},
\label{tau}
\end{align}
where $C_P$ is given in \eqref{eq2}, $\widehat{C}$ is defined in \eqref{sobolev} }\extra{ and $B_r$ is a ball centered at $\textbf{0}$ of radius $r>0$ sufficiently large such that $x-\Omega\subset B_r$ for any $x\in \Omega$} (see also \cite{g} for this observation on $B_r$).
\ry{ Now observe that, since, by \eqref{F}, there exists a positive constant $C_F>0$ such that, for $\delta$ sufficiently small,  $$0<\frac{1}{\Fsec(1-2\delta)}\leq C_F\delta \quad\text{ and }\quad0<\frac{1}{F^\prime(1-2\delta)}\leq \frac{C_F}{\vert\ln(\delta)\vert},$$ we have 
 \begin{align*}
    \dfrac{\frac{8\delta^{2}}{\widetilde{\tau}}}{\frac{\Vert\nabla J\Vert_{{L^1(B_r)}}^2}{2\Fsec(1-2\delta)}}&=\dfrac{16\delta^2\Fsec(1-2\delta)}{\Vert\nabla J\Vert_{{L^1(B_r)}}^2}\dfrac{3C(\tau)\Vert\nabla J\Vert_{{L^1(B_r)}}^{5}\widehat{C}^{\frac{3}{2}}\left(1+C_P^2\right)^{\frac{3}{2}}}{2^{-20}\delta^5\left(\Fsec(1-2\delta)\right)^{4}F^\prime(1-2\delta)}\\&=\dfrac{3C(\tau)\Vert\nabla J\Vert_{{L^1(B_r)}}^{3}\widehat{C}^{\frac{3}{2}}\left(1+C_P^2\right)^{\frac{3}{2}}}{2^{-24}\delta^3\left(\Fsec(1-2\delta)\right)^{3}F^\prime(1-2\delta)}\leq \dfrac{\widetilde{C}}{\vert \ln(\delta)\vert}\to 0 \quad\text{as }\delta\to 0^+,
\end{align*}
where $\widetilde{C}:=\dfrac{3C(\tau)\Vert\nabla J\Vert_{{L^1(B_r)}}^{3}\widehat{C}^{\frac{3}{2}}\left(1+C_P^2\right)^{\frac{3}{2}}C_F^4}{2^{-24}}>0,$ so that 
\begin{align*}
    \dfrac{\frac{8\delta^{2}}{\widetilde{\tau}}}{\frac{\Vert\nabla J\Vert_{{L^1(B_r)}}^2}{2\Fsec(1-2\delta)}}=O\left(\dfrac{1}{\vert\ln(\delta)\vert}\right).
\end{align*}
 }
\ry{This} means that we can find a sufficiently small $\delta>0$ so that 
\begin{align}
\max\left\{\frac{\Vert\nabla J\Vert_{\extra{L^1(B_r)}}^2}{2\Fsec(1-2\delta)},\frac{8\delta^2}{\widetilde{\tau}}\right\}=\frac{\Vert\nabla J\Vert_{\extra{L^1(B_r)}}^2}{2\Fsec(1-2\delta
)}.\label{delt}
\end{align}
Choose now  $T>0$ such that $T-3\widetilde{\tau}\geq \frac{\tau}{2}$ (for example, one can start with $T=3\widetilde{\tau}+\frac{\tau}{2}$). Up to reducing the size of $\delta$, and thus of $\widetilde{\tau}$, we can find $\widetilde{\tau}$ such that 
\begin{equation}
2\widetilde{\tau}+\frac{\tau}{2}\leq \tau.
\label{tt}
\end{equation}

Let us then fix $\delta>0$ (and thus $\widetilde{\tau}>0$) so that also \eqref{delt} and  \eqref{tt} hold. Notice that the choice of $\delta$ and $\widetilde{\tau}$ \textit{does not} depend on the specific $T$, but clearly depends on $\tau$.\\
We now define the sequence
\begin{align}
k_n=1-\delta-\frac{\delta}{2^n}, \quad \forall n\geq 0,
\label{kn}
\end{align}
where 
\begin{align}
1-2\delta< k_n<k_{n+1}<1-\delta,\qquad \forall n\geq 1,\qquad k_n\to 1-\delta\qquad \text{as }n\to \infty,
\label{kn1}
\end{align}
and the sequence of times
\begin{align}
\begin{cases}
t_{-1}=T-3\widetilde{\tau},\\
t_n=t_{n-1}+\frac{\widetilde{\tau}}{2^n},\qquad n\geq 0,
\end{cases}
\end{align}
 which satisfies
$$
t_{-1}<t_n<t_{n+1}< T-\widetilde{\tau},\qquad \forall n\geq 0.
$$
We now introduce a cutoff function $\eta_n\in C^1(\R)$ by setting 
\begin{align}
\eta_n(t):=\begin{cases}
0,\quad t\leq t_{n-1},\\
1,\quad t\geq t_{n},
\end{cases}\text{ and }\quad \vert \eta^\prime_n(t)\vert\leq \frac{2^{n+1}}{\widetilde{\tau}},
\label{cutoff}
\end{align}
on account of the above definition of the sequence $\{t_n\}_n$.
Recalling \eqref{phik}, we then set \rx{$\rho=k_n$},
\begin{align}
\phi_n(x,t):=(\phi-k_n)^+,
\label{phik0}
\end{align}
and, for any $n\geq 0$, we introduce the interval $I_n=[t_{n-1},T]$ and the set
$$
A_n(t):=\{x\in \Omega: \phi(x,t)-k_n\geq 0\},\quad \forall t\in I_n.
$$
Clearly, we have
$$
I_{n+1}\subseteq I_n,\qquad \forall n\geq 0,$$
$$A_{n+1}(t)\subseteq A_n(t),\qquad \forall n\geq 0,\qquad \forall t\in I_{n+1}.
$$
In conclusion, we set
$$
y_n=\int_{I_n}\int_{A_n(s)}1dxds,\qquad \forall n\geq0.
$$
 Now, for any $n\geq 0$, we consider the test function $v=\phi_n\eta_n^2$, and integrate over $[t_{n-1},t]$, $t_n\leq t\leq T$. Then we have 
\begin{align}
&\nonumber\int_{t_{n-1}}^t<\partial_t\phi,\phi_n\eta_n^2>ds+\int_{t_{n-1}}^t\int_{A_n(s)}F^{\prime\prime}(\phi)\nabla\phi\cdot \nabla\phi_n \eta_n^2 dx ds\\&=\int_{t_{n-1}}^t\int_{A_n(s)}\eta_n^2(\nabla J\ast \phi)\cdot \nabla \phi_ndx ds,
\label{phin}
\end{align}
since $\nabla F^\prime(\phi(t))=\Fsec(\phi)\nabla\phi(t)$, for almost every $x\in\Omega$ and for any $t\geq \frac{\tau}{2}$, which can be proven, e.g., by a truncation argument as in \cite[Lemma 3.2]{Wu}, applied for any $t\geq \frac{\tau}{2}$. Indeed, as in \cite[(3.5)]{Wu}, we obtain $\nabla F^\prime(\phi(t))=\Fsec(\phi)\nabla\phi(t)$ in the sense of distribution and thus, being $\nabla F^\prime(\phi)\in L^\infty(\frac{\tau}{2},\infty;\textbf{H})$, we immediately infer that the equality holds also almost everywhere in $\Omega$, for any $t\geq \frac{\tau}{2}$.
Now, as in \cite{GGG}, \extra{for $\delta$ sufficiently small} we obtain
\begin{align}
\int_{t_{n-1}}^t\eta_n^2\int_{A_n(s)}F^{\prime\prime}(\phi)\nabla\phi\cdot \nabla\phi_n dx ds\geq \Fsec(1-2\delta)\int_{t_{n-1}}^t\eta_n^2\Vert \nabla\phi_n\Vert^2ds,
\label{Fss}
\end{align}
and, for the right-hand side of \eqref{phin}, recalling that $\vert\phi\vert< 1$ a.e. in $\Omega\times (0,+\infty)$, we find
\begin{align}
&\nonumber\int_{t_{n-1}}^t\int_{A_n(s)}(\nabla J\ast \phi)\cdot \nabla \phi_n\eta_n^2dx ds\\&\nonumber\leq \frac{1}{2}\Fsec(1-2\delta)\int_{t_{n-1}}^t\eta_n^2\Vert \nabla\phi_n\Vert^2ds+\frac{1}{2\Fsec(1-2\delta)}\int_{t_{n-1}}^t\int_{A_n(s)}\eta_n^2\vert \nabla J\ast \phi\vert^2 dxds\\&\nonumber
\leq 
\frac{1}{2}\Fsec(1-2\delta)\int_{t_{n-1}}^t\eta_n^2\Vert \nabla\phi_n\Vert^2ds+\frac{1}{2\Fsec(1-2\delta)}\int_{t_{n-1}}^t\Vert\nabla J\ast \phi\Vert_{L^\infty(\Omega)}^2\int_{A_n(s)}dxds\\&\nonumber\leq 
\frac{1}{2}\Fsec(1-2\delta)\int_{t_{n-1}}^t\eta_n^2\Vert \nabla\phi_n\Vert^2ds+\frac{\Vert\nabla J\Vert_{\extra{L^1(B_r)}}^2}{2\Fsec(1-2\delta)}\int_{t_{n-1}}^t\int_{A_n(s)}dxds\\&\leq \frac{1}{2}\Fsec(1-2\delta)\int_{t_{n-1}}^t \eta_n^2\Vert \nabla\phi_n\Vert^2ds+\frac{\Vert\nabla J\Vert_{\extra{L^1(B_r)}}^2}{2\Fsec(1-2\delta)}y_n,
\label{J0}
\end{align}
where we have applied (see, e.g., \cite[Thm. 4.33]{Brezis}): 
\begin{align}
\Vert\nabla J\ast \phi\Vert_{L^\infty(\Omega)}\leq \Vert\nabla J\Vert_{\extra{L^1(B_r)}}\Vert\phi\Vert_{L^\infty(\Omega)}\leq \Vert\nabla J\Vert_{\extra{L^1(B_r)}}.
\label{J}
\end{align}
Moreover, we have 
\begin{align}
\int_{t_{n-1}}^t<\partial_t\phi,\phi_n\eta_n^2>ds=\frac{1}{2}\Vert\phi_n(t)\Vert^2-\int_{t_{n-1}}^t\Vert\phi_n(s)\Vert^2\eta_n\partial_t\eta_nds.
\label{eq1z}
\end{align}
Note that, since $ \vert\phi\vert<1\text{ a.e. in }\Omega$, for any $t\geq \frac{\tau}{2}$, 
\begin{align}
0\leq \phi_n\leq 2\delta\quad\text{ a.e. in }\Omega,\quad \forall t\geq \frac{\tau}{2}.
\label{essential}
\end{align}
Then, by the above inequality,
\begin{align}
\nonumber
&\int_{t_{n-1}}^t\Vert\phi_n(s)\Vert^2\eta_n\partial_t\eta_nds=\int_{t_{n-1}}^t\int_\Omega \phi_n^2(s)\eta_n\partial_t\eta_ndxds=\int_{t_{n-1}}^t\int_{A_n(s)} \phi_n^2(s)\eta_n\partial_t\eta_ndxds\\&
\leq 
\int_{t_{n-1}}^t\int_{A_n(s)} (2\delta)^2\frac{2^{n+1}}{\widetilde{\tau}}dxds\leq \frac{2^{n+3}\delta^2}{\widetilde{\tau}}y_n.
\label{del}
\end{align}
\rx{Plugging \eqref{Fss}, \eqref{J0}, \eqref{eq1z} and \eqref{del} into \eqref{phin}, we find}
\begin{align*}
&\frac{1}{2}\Vert\phi_n(t)\Vert^2+\frac{1}{2}\Fsec(1-2\delta)\int_{t_{n-1}}^t \eta_n^2\Vert \nabla\phi_n(s)\Vert^2ds\\&\leq2^{n+1}\max\left\{\frac{\Vert\nabla J\Vert_{\extra{L^1(B_r)}}^2}{2\Fsec(1-2\delta)},\frac{8\delta^2}{\widetilde{\tau}}\right\}y_n,
\end{align*}
for any $t\in[t_n,T]$. Thanks to the choice of $\delta$ and $\widetilde{\tau}$, we recall \eqref{delt}, implying
\begin{align}
\max_{t\in I_{n+1}}\Vert\phi_n(t)\Vert^2\leq X_n,\qquad  \Fsec(1-2\delta)\int_{I_{n+1}}\Vert \nabla\phi_n\Vert^2ds \leq X_n,
\label{est}
\end{align}
where 
$$
X_n:= 2^{n+1}\frac{\Vert\nabla J\Vert_{\extra{L^1(B_r)}}^2}{\Fsec(1-2\delta)}y_n.
$$
On the other hand, for any $t\in I_{n+1}$ and for almost any $x\in A_{n+1}(t)$, we get
\begin{align*}
&\phi_n(x,t)=\phi(x,t)-\left[1-\delta-\frac{\delta}{2^n}\right]\\&=
\underbrace{\phi(x,t)-\left[1-\delta-\frac{\delta}{2^{n+1}}\right]}_{\phi_{n+1}(x,t)\geq 0}+\delta\left[\frac{1}{2^{n}}-\frac{1}{2^{n+1}}\right]\geq \frac{\delta}{2^{n+1}},
\end{align*}
which implies
\begin{align*}
\int_{I_{n+1}}\int_{\Omega}\vert\phi_n\vert^3dxds\geq \int_{I_{n+1}}\int_{A_{n+1}(s)}\vert\phi_n\vert^3dxds\geq \left(\frac{\delta}{2^{n+1}}\right)^3\int_{I_{n+1}}\int_{A_{n+1}(s)}dxds=\left(\frac{\delta}{2^{n+1}}\right)^3y_{n+1}.
\end{align*}
Then we have 
\begin{align}
&\nonumber\left(\frac{\delta}{2^{n+1}}\right)^3y_{n+1}\leq \int_{I_{n+1}}\int_{\Omega}\vert\phi_n\vert^3dxds\\&=\int_{I_{n+1}}\int_{A_n(s)}\vert\phi_n\vert^3dxds \leq \left(\int_{I_{n+1}}\int_{\Omega}\vert\phi_n\vert^{\frac{10}{3}}dxds\right)^{\frac{9}{10}}\left(\int_{I_{n+1}}\int_{A_n(s)}dxds\right)^{\frac{1}{10}}.
\label{est2}
\end{align}
Notice that, thanks to \eqref{sobolev} and \eqref{eq2} (which holds thanks to \eqref{kn1}), we get 
\begin{align*}
&\int_{I_{n+1}}\int_{\Omega}\vert\phi_n\vert^{\frac{10}{3}}dxds\leq \widehat{C} \int_{I_{n+1}}\Vert\phi_n\Vert_{V}^{2}\Vert\phi_n\Vert^{\frac{4}{3}}ds\leq \widehat{C}\int_{I_{n+1}} \left(\Vert\phi_n\Vert^2+\Vert\nabla\phi_n\Vert^2\right)\Vert\phi_n\Vert^{\frac{4}{3}}ds\\&\leq \widehat{C}\left(1+C_P^2\right)\int_{I_{n+1}}\Vert\nabla\phi_n\Vert^2\Vert\phi_n\Vert^{\frac{4}{3}}ds,
\end{align*}
where we have chosen an equivalent norm on $V$.
Observe now that, by \eqref{est}, 
\begin{align*}
&\int_{I_{n+1}}\int_{\Omega}\vert\phi_n\vert^{\frac{10}{3}}dxds\leq \widehat{C}\left(1+C_P^2\right)\int_{I_{n+1}}\Vert\nabla\phi_n\Vert^2\Vert\phi_n\Vert^{\frac{4}{3}}ds\leq \widehat{C}\left(1+C_P^2\right)\max_{t\in I_{n+1}}\Vert\phi_n(t)\Vert^{\frac{4}{3}}\int_{I_{n+1}}\Vert\nabla\phi_n\Vert^2ds\\&
\leq \frac{\widehat{C}(1+C_P^2)}{\Fsec(1-2\delta)}X_n^{\frac{2}{3}}\Fsec(1-2\delta)\int_{I_{n+1}}\Vert\nabla\phi_n\Vert^2ds\leq \frac{\widehat{C}(1+C_P^2)}{\Fsec(1-2\delta)}X_n^{\frac{5}{3}}\leq \frac{2^{\frac{5n}{3}+\frac{5}{3}}\Vert\nabla J\Vert_{\extra{L^1(B_r)}}^{\frac{10}{3}}\widehat{C}\left(1+C_P^2\right)}{(\Fsec(1-2\delta))^{\frac{8}{3}}}y_n^{\frac{5}{3}}.
\end{align*}
 
Coming back to \eqref{est2}, we immediately infer
\begin{align}
&\nonumber\left(\frac{\delta}{2^{n+1}}\right)^3y_{n+1}\leq \left(\int_{I_{n+1}}\int_{\Omega}\vert\phi_n\vert^{\frac{10}{3}}dxds\right)^{\frac{9}{10}}\left(\int_{I_{n+1}}\int_{A_n(s)}dxds\right)^{\frac{1}{10}}\\&\leq \frac{2^{{\frac{3}{2}n+\frac{3}{2}}}\Vert\nabla J\Vert_{\extra{L^1(B_r)}}^{3}\widehat{C}^{\frac{9}{10}}\left(1+C_P^2\right)^{\frac{9}{10}}}{\left(\Fsec(1-2\delta)\right)^{\frac{12}{5}}}y_n^{\frac{3}{2}}y_n^{\frac{1}{10}}=\frac{2^{{\frac{3}{2}n+\frac{3}{2}}}\Vert\nabla J\Vert_{\extra{L^1(B_r)}}^{3}\widehat{C}^{\frac{9}{10}}\left(1+C_P^2\right)^{\frac{9}{10}}}{\left(\Fsec(1-2\delta)\right)^{\frac{12}{5}}}y_n^{\frac{8}{5}}.
\label{est3}
\end{align}
In conclusion, we end up with 
\begin{align}
y_{n+1}\leq \frac{2^{{\frac{9}{2}n+\frac{9}{2}}}\Vert\nabla J\Vert_{\extra{L^1(B_r)}}^{3}\widehat{C}^{\frac{9}{10}}\left(1+C_P^2\right)^{\frac{9}{10}}}{\delta^3\left(\Fsec(1-2\delta)\right)^{\frac{12}{5}}}y_n^{\frac{8}{5}},\qquad \forall n\geq 0.
\label{last0}
\end{align}
Thus we can apply Lemma \ref{conv}. In particular, we have $b=2^\frac{9}{2}>1$, $C=\frac{2^{\frac{9}{2}}\Vert\nabla J\Vert_{\extra{L^1(B_r)}}^{3}\widehat{C}^{\frac{9}{10}}\left(1+C_P^2\right)^{\frac{9}{10}}}{\delta^3\left(\Fsec(1-2\delta)\right)^{\frac{12}{5}}}>0$, $\varepsilon=\frac{3}{5}$, to get that ${y}_n\to 0$, as long as 
$$
{y}_0\leq C^{-\frac{5}{3}}b^{-\frac{25}{9}},
$$
i.e.,
\begin{align}
y_0\leq \dfrac{2^{-20}\delta^5\left(\Fsec(1-2\delta)\right)^{4}}{\Vert\nabla J\Vert_{\extra{L^1(B_r)}}^{5}\widehat{C}^{\frac{3}{2}}\left(1+C_P^2\right)^{\frac{3}{2}}}.
\label{last}
\end{align}
We are left with a last estimate: thanks to \eqref{Fprime}, we know that $\Vert F^\prime(\phi)\Vert_{L^\infty(\frac{\tau}{2},\infty;{L^1(\Omega)})}\leq C(\tau)$  and $F^\prime$ is monotone in a neighborhood of $+1$, so that we infer
\begin{align*}
&y_0=\int_{I_0}\int_{A_0(s)}1dxds\leq\int_{I_0}\int_{\{x\in\Omega:\ \phi(x,s) \geq 1-2\delta\}}1dxds\\&\leq \int_{I_0}\int_{A_0(s)}\frac{\vert F^\prime(\phi)\vert}{F^\prime(1-2\delta)} dxds\leq \frac{3C(\tau)\widetilde{\tau}}{{F^\prime(1-2\delta)}}.
\end{align*}
Therefore, if we ensure that 
$$
\frac{3C(\tau)\widetilde{\tau}}{{F^\prime(1-2\delta)}}\leq \dfrac{2^{-20}\delta^5\left(\Fsec(1-2\delta)\right)^{4}}{\Vert\nabla J\Vert_{\extra{L^1(B_r)}}^{5}\widehat{C}^{\frac{3}{2}}\left(1+C_P^2\right)^{\frac{3}{2}}},
$$
then \eqref{last} holds. Having fixed $\widetilde{\tau}$ \rx{in \eqref{tau}} such that
\begin{align}
\widetilde{\tau}= \dfrac{2^{-20}\delta^5\left(\Fsec(1-2\delta)\right)^{4}F^\prime(1-2\delta)}{3C(\tau)\Vert\nabla J\Vert_{\extra{L^1(B_r)}}^{5}\widehat{C}^{\frac{3}{2}}\left(1+C_P^2\right)^{\frac{3}{2}}},
\label{taubis}
\end{align}
we obtain the result. Notice that $\delta$ is fixed, so $\widetilde{\tau}>0$ is not infinitesimal, but it depends on  $\phi_0$ in a nontrivial way (thus not only on the initial energy) through $C_P$.

In the end, passing to the limit in $y_n$ as $n\to\infty$, we have obtained that 
$$
\Vert(\phi-(1-\delta))^+\Vert_{L^\infty(\Omega\times({T}-\widetilde{\tau},{T}))}=0,
$$
since, as $n\to\infty$, 
\begin{align*}
    y_n\to \left\vert\left\{(x,t)\in \Omega\times[T-\widetilde{\tau},T]: \phi(x,t)\geq 1-\delta\right\}\right\vert=0.
\end{align*}
 We now repeat exactly the same argument for the case $(\phi-(-1+\delta))^-$ (using $\phi_n(t)=(\phi(t)+k_n)^-$).
Notice that also for this second case we have the same constant $C_P$ (see \eqref{eq2}). Moreover, the argument is exactly the same due to assumption \eqref{F2}, which implies that $\frac{1}{\Fsec(-1+2\delta)}=O(\delta)$ and $\frac{1}{\vert F^\prime(-1+2\delta)\vert }=O\left(\frac{1}{\vert\ln(\delta)\vert}\right)$, for $\delta$ sufficiently small. We can then choose the minima between the $\delta$ and $\widetilde{\tau}$ obtained in the two cases, to get in the end that there exists a couple $\delta>0,\widetilde{\tau}>0$ such that
 \begin{align}
-1+\delta\leq \phi(x,t)\leq 1-\delta,\quad\text{a.e. in }\Omega\times ({T}-\widetilde{\tau},{T}).
\label{end}
\end{align}
Finally, notice that, due to the choice of $T$, we have $T-\widetilde{\tau}=2\widetilde{\tau}+\frac{\tau}{2}\leq \tau$, therefore we can repeat the same procedure on the interval $({T},T+\widetilde{\tau})$ (this means that the new starting time will be $t_{-1}={T}-2\widetilde{\tau}\geq \frac{\tau}{2}$) and so on, reaching eventually the entire interval $[\tau,+\infty)$. Clearly $\delta$ and $\widetilde{\tau}$ are always the same, since the constant $C_P$ is uniform over the entire interval $[\tau,+\infty)$ and the time horizon $T$ does not enter in any of the estimates. The proof is thus concluded.
\begin{remark}
	We point out that the same proof holds for the case of convective nonlocal Cahn-Hilliard equation:
	 \begin{equation}
	\label{nonloc3d_conv}
	\begin{cases}
	\partial_t\phi+\textbf{u}\cdot \nabla\phi-\Delta\mu=0\quad\text{in }\Omega\times(0,T),\\
	\mu=F^\prime(\phi)-J\ast \phi\quad\text{in }\Omega\times(0,T),\\
	\partial_\textbf{n}\mu=0\quad\text{ on }\partial\Omega\times(0,T),\\
	\phi(\cdot,0)=\phi_0\quad \text{ in }\Omega,
	\end{cases}
	\end{equation}
	where $\textbf{u}$ is a sufficiently regular divergence free  vector field, such that $\textbf{u}\cdot \textbf{n}=0$ on $\partial\Omega\times(0,T)$. Indeed, Theorem \ref{known} can be mostly extended also to this case (see, e.g., \cite[Sec.6]{GGG0}, in which a related system, the nonlocal Cahn-Hilliard-Navier-Stokes system, is analyzed). Moreover, in the proof of Theorem \ref{sep} the term $\textbf{u}\cdot \nabla\phi$ does not appear, since in \eqref{phin} we should get an additional $(\textbf{u}\cdot\nabla\phi, \phi_n\eta_n^2)$, which is zero thanks to the assumptions on $\textbf{u}$. Therefore, the separation property could a priori be obtained also in the couplings of the nonlocal Cahn-Hilliard equation with some hydrodynamic models, like Navier-Stokes equations (see, e.g., \cite{Abelsterasawa} or \cite[Sec.6]{GGG0} for some examples of such models).
	\label{extension}
\end{remark}
\begin{remark}
One might think that the proof of Theorem \ref{sep} could be adapted also to the conserved Allen-Cahn equation
    \begin{equation}
	\label{nonloc3d_conv1}
	\begin{cases}
	\partial_t\phi+\mu-\overline\mu=0\quad\text{in }\Omega\times(0,T),\\
	\mu=\Psi^\prime(\phi)-\Delta \phi\quad\text{in }\Omega\times(0,T),\\
	\partial_\textbf{n}\phi=0\quad\text{ on }\partial\Omega\times(0,T),\\
	\phi(\cdot,0)=\phi_0\quad \text{ in }\Omega,
	\end{cases}
	\end{equation}
	where $\Psi$ is defined in \eqref{potential_loc}.
	\extra{Indeed, this has been obtained in the recent \cite{multiAC} in the case of multi-component conserved Allen-Cahn equation in two and three dimensions, and it is valid also for \eqref{nonloc3d_conv1}. In the proof one loses the term $$\int_{t_{n-1}}^t\int_{A_n(s)}F^{\prime\prime}(\phi)\nabla\phi\cdot \nabla\phi_n \eta_n^2 dx ds,$$ which is substituted by $\int_{t_{n-1}}^t\int_{A_n(s)}F^{\prime}(\phi)\phi_n \eta_n^2 dx ds\geq F'(1-2\delta)\int_{t_{n-1}}^t\int_{A_n(s)}\phi_n \eta_n^2 dx ds$: the presence of the first derivative of $F$ instead of the second derivative, since $F'(1-2\delta)\to+\infty$ as $\delta\to 0^+$, is still enough to carry out the De Giorgi's iteration scheme, by heavily exploiting estimate \eqref{essential}. We also mention the fact that in two-dimensional bounded domains the instantaneous strict separation property for \eqref{nonloc3d_conv1} was proven before in \cite{GGW}, by a completely different argument.} 
\end{remark}
\ry{\begin{remark}
Assumption \ref{h4} shows that the strict separation property also holds for more general and singular potentials $F$ than the logarithmic one \eqref{potential_loc}. Furthermore, by slightly adapting the proof of  Theorem \ref{sep}, one can show that the same property also holds for more general double well potentials than $F$. For instance, one could deal with a chemical potential $\mu=\Psi^\prime(\phi)+(J\ast 1)\phi-J\ast \phi$, with $\Psi$ defined in \eqref{potential_loc} and obtain an analogous result. Notice that in this new setting the nonlocal term $J\ast \phi$ is related to diffusion effects (see \cite{GGG} and the references therein for more details). Also, in the case of non-constant mobility $M(\phi)$, the proof should work well as long as it is nondegenerate (i.e., bounded below by a strictly positive constant) and the existence of strong solutions is given. In conclusion, another possible extension could be in the case of dynamic boundary conditions (see, e.g., \cite{PaS}): first one needs to assess the existence of strong solutions and the instantaneous regularization of weak solutions, and then apply a De Giorgi's iteration scheme, which seems harder due to the presence of boundary terms which have to be carefully handled.  
\end{remark}}
\subsection{Proof of Corollary \ref{sep2}}
Observe that, due to Remark \ref{proof}, we only need to prove that the unique solution $\phi$ departing from $\phi_0$ is strictly separated from the pure phases in a neighborhood of the initial time. To this aim we perform again a De Giorgi's iteration scheme, in this case without the use of a cutoff function of the form \eqref{cutoff}. Indeed, the necessity of the cutoff function is merely to eliminate the presence of the initial datum in the estimates, but in our case, up to choosing $\delta\leq \frac{\delta_0}{2}$, this problem does not appear any more, as we shall see. Again the Step 1 of the proof of Theorem \ref{sep} is still valid, and we adopt the same notation. Clearly, thanks to Remark \ref{strong}, we can choose $\tau=0$, \rx{so that again}
		\begin{equation}
\Vert \phi_{\rx{\rho}}(t)\Vert\leq C_{P}\Vert\nabla \phi_{\rx{\rho}}(t)\Vert, \qquad 	\Vert \widetilde{\phi}_{\rx{\rho}}(t)\Vert\leq C_{P}\Vert\nabla \widetilde{\phi}_{\rx{\rho}}(t)\Vert,\quad \text{ for almost any }t\geq 0\text{ and for any }\rx{\rho\in[1-2\widehat{\delta},1].}
\label{eq3}
\end{equation}
We then start from Step 2. Let us fix $\delta$ sufficiently small such that $\delta\leq \min\left\{\widehat{\delta},\frac{\delta_0}{2}\right\}$, so that \eqref{eq3} holds for any $\rx{\rho}\in[ 1-2\delta,1]$. Set then $\widetilde{\tau}>0$ such that \eqref{tau2} below holds. 
As in Theorem \ref{sep}, we define the same sequence \eqref{kn}, but we do not need to consider any sequence of times, since we will always use the same, fixed, interval $I:=[0,\widetilde{\tau}]$.
Then we define again
\begin{align}
\phi_n(x,t):=(\phi-k_n)^+,
\label{phiks}
\end{align}
and, for any $n\geq 0$, we introduce the set 
$$
A_n(t):=\{x\in \Omega: \phi(x,t)-k_n\geq 0\},\quad \forall t\in I,
$$
so that
$$A_{n+1}(t)\subseteq A_n(t),\qquad \forall n\geq 0,\qquad \forall t\in I.
$$
We thus set
$$
y_n=\int_{I}\int_{A_n(s)}1dxds,\qquad \forall n\geq0.
$$
Now, for any $n\geq 0$ we consider the test function $w=\phi_n$, and integrate over $[0,t]$, $ t\leq \widetilde{\tau}$. Then we have, as in Theorem \ref{sep},
\begin{align*}
\frac{1}{2}\Vert\phi_n(t)\Vert^2+\int_{0}^t\int_{A_n(s)}F^{\prime\prime}(\phi)\nabla\phi\cdot \nabla\phi_n dx ds=\int_{0}^t\int_{A_n(s)}(\nabla J\ast \phi)\cdot \nabla \phi_ndx ds+\frac{1}{2}\Vert\phi_n(0)\Vert^2.
\end{align*}
Note that, due to the choice of $\delta\leq \frac{\delta_0}{2}$, thanks to the strict separation of the initial datum, we immediately infer that $\Vert\phi_n(0)\Vert=0$ for any $n\geq 0$.
Following the same arguments as in the proof of Theorem \ref{sep}, we obtain 
\begin{align*}
&\frac{1}{2}\Vert\phi_n(t)\Vert^2+\frac{1}{2}\Fsec(1-2\delta)\int_{0}^t \Vert \nabla\phi_n(s)\Vert^2ds\leq\frac{\Vert\nabla J\Vert_{\extra{L^1(B_r)}}^2}{2\Fsec(1-2\delta)}y_n,
\end{align*}
for any $t\in[0,\widetilde{\tau}]$. Observe that we do not see the presence of the term related to $\frac{1}{\widetilde{\tau}}$ (estimated in \eqref{del}), since it is a consequence of the use of the cutoff function \eqref{cutoff}. This implies
\begin{align}
\max_{t\in I}\Vert\phi_n(t)\Vert^2\leq Z_n,\qquad  \Fsec(1-2\delta)\int_{I}\Vert \nabla\phi_n\Vert^2ds \leq Z_n,
\label{est4}
\end{align}
where 
$$
\ry{Z_n:=\frac{\Vert\nabla J\Vert_{{L^1(B_r)}}^2}{\Fsec(1-2\delta)}y_n.}
$$
Observe that, for any $t\in I$ and for almost any $x\in A_{n+1}(t)$, we get
\begin{align*}
&\phi_n(x,t)=
\underbrace{\phi(x,t)-\left[1-\delta-\frac{\delta}{2^{n+1}}\right]}_{\phi_{n+1}(x,t)\geq 0}+\delta\left[\frac{1}{2^{n}}-\frac{1}{2^{n+1}}\right]\geq \frac{\delta}{2^{n+1}},
\end{align*}
which implies
\begin{align*}
\int_{I}\int_{\Omega}\vert\phi_n\vert^3dxds\geq \int_{I}\int_{A_{n+1}(s)}\vert\phi_n\vert^3dxds\geq \left(\frac{\delta}{2^{n+1}}\right)^3\int_{I}\int_{A_{n+1}(s)}dxds=\left(\frac{\delta}{2^{n+1}}\right)^3y_{n+1}.
\end{align*}
Then we have, as in \eqref{est2},
\begin{align}
&\left(\frac{\delta}{2^{n+1}}\right)^3y_{n+1}\leq   \left(\int_{I}\int_{\Omega}\vert\phi_n\vert^{\frac{10}{3}}dxds\right)^{\frac{9}{10}}\left(\int_{I}\int_{A_n(s)}dxds\right)^{\frac{1}{10}}.
\label{est1}
\end{align}
Again thanks to \eqref{sobolev} and \eqref{eq3}, we have 
\begin{align*}
&\int_{I}\int_{\Omega}\vert\phi_n\vert^{\frac{10}{3}}dxds\leq \widehat{C}\left(1+C_P^2\right)\int_{I}\Vert\nabla\phi_n\Vert^2\Vert\phi_n\Vert^{\frac{4}{3}}ds,
\end{align*}
so that, by \eqref{est4}, 
\begin{align*}
&\int_{I}\int_{\Omega}\vert\phi_n\vert^{\frac{10}{3}}dxds\leq \widehat{C}\left(1+C_P^2\right)\int_{I}\Vert\nabla\phi_n\Vert^2\Vert\phi_n\Vert^{\frac{4}{3}}ds\leq \widehat{C}\left(1+C_P^2\right)\max_{t\in I}\Vert\phi_n\Vert^{\frac{4}{3}}\int_{I}\Vert\nabla\phi_n\Vert^2ds\\&
\leq \frac{\widehat{C}(1+C_P^2)}{\Fsec(1-2\delta)}Z_n^{\frac{2}{3}}\Fsec(1-2\delta)\int_{I}\Vert\nabla\phi_n\Vert^2ds\leq \frac{\widehat{C}(1+C_P^2)}{\Fsec(1-2\delta)}Z_n^{\frac{5}{3}}\leq \ry{\frac{\Vert\nabla J\Vert_{{L^1(B_r)}}^{\frac{10}{3}}\widehat{C}\left(1+C_P^2\right)}{(\Fsec(1-2\delta))^{\frac{8}{3}}}y_n^{\frac{5}{3}}}.
\end{align*}
Therefore, we immediately infer from \eqref{est1}
\begin{align}
&\nonumber\left(\frac{\delta}{2^{n+1}}\right)^3y_{n+1}\leq \left(\int_{I}\int_{\Omega}\vert\phi_n\vert^{\frac{10}{3}}dxds\right)^{\frac{9}{10}}\left(\int_{I}\int_{A_n(s)}dxds\right)^{\frac{1}{10}}\\&\leq \ry{\frac{\Vert\nabla J\Vert_{{L^1(B_r)}}^{3}\widehat{C}^{\frac{9}{10}}\left(1+C_P^2\right)^{\frac{9}{10}}}{\left(\Fsec(1-2\delta)\right)^{\frac{12}{5}}}y_n^{\frac{3}{2}}y_n^{\frac{1}{10}}=\frac{\Vert\nabla J\Vert_{{L^1(B_r)}}^{3}\widehat{C}^{\frac{9}{10}}\left(1+C_P^2\right)^{\frac{9}{10}}}{\left(\Fsec(1-2\delta)\right)^{\frac{12}{5}}}y_n^{\frac{8}{5}}.}
\label{est3bis}
\end{align}
In conclusion, we end up with 
\begin{align*}
y_{n+1}\leq \frac{2^{{3n+\ry{3}}}\Vert\nabla J\Vert_{\extra{L^1(B_r)}}^{3}\widehat{C}^{\frac{9}{10}}\left(1+C_P^2\right)^{\frac{9}{10}}}{\delta^3\left(\Fsec(1-2\delta)\right)^{\frac{12}{5}}}y_n^{\frac{8}{5}},\qquad \forall n\geq 0,
\end{align*}
and we can apply Lemma \ref{conv}. In particular, we have $b=2^3>1$, $C=\frac{2^{\ry{3}}\Vert\nabla J\Vert_{\extra{L^1(B_r)}}^{3}\widehat{C}^{\frac{9}{10}}\left(1+C_P^2\right)^{\frac{9}{10}}}{\delta^3\left(\Fsec(1-2\delta)\right)^{\frac{12}{5}}}>0$, $\varepsilon=\frac{3}{5}$, to get that ${y}_n\to 0$, as long as 
$$
{y}_0\leq C^{-\frac{5}{3}}b^{-\frac{25}{9}},
$$
i.e.,
\begin{align}
y_0\leq \dfrac{2^{-\ry{\frac{40}{3}}}\delta^5\left(\Fsec(1-2\delta)\right)^{4}}{\Vert\nabla J\Vert_{\extra{L^1(B_r)}}^{5}\widehat{C}^{\frac{3}{2}}\left(1+C_P^2\right)^{\frac{3}{2}}}.
\label{lastbis}
\end{align}
In conclusion, since we know by \eqref{Fprime} and Remark \ref{strong} that $\Vert F^\prime(\phi)\Vert_{L^\infty(0,\infty;{L^1(\Omega)})}\leq C$, we infer
\begin{align*}
&y_0=\int_{I}\int_{A_0(s)}1dxds\leq \int_{I}\int_{A_0(s)}\frac{\vert F^\prime(\phi)\vert}{F^\prime(1-2\delta)} dxds\leq \frac{C\widetilde{\tau}}{{F^\prime(1-2\delta)}}.
\end{align*}
Having fixed $\widetilde{\tau}$ so that
\begin{align}
\widetilde{\tau}= \dfrac{2^{-\ry{\frac{40}{3}}}\delta^5\left(\Fsec(1-2\delta)\right)^{4}F^\prime(1-2\delta)}{C\Vert\nabla J\Vert_{\extra{L^1(B_r)}}^{5}\widehat{C}^{\frac{3}{2}}\left(1+C_P^2\right)^{\frac{3}{2}}},
\label{tau2}
\end{align}
we have
$$
\frac{C\widetilde{\tau}}{{F^\prime(1-2\delta)}}\leq \dfrac{2^{-\ry{\frac{40}{3}}}\delta^5\left(\Fsec(1-2\delta)\right)^{4}}{\Vert\nabla J\Vert_{\extra{L^1(B_r)}}^{5}\widehat{C}^{\frac{3}{2}}\left(1+C_P^2\right)^{\frac{3}{2}}},
$$
so that \eqref{lastbis} holds. 
In the end, passing to the limit in $y_n$ as $n\to\infty$, we have obtained that 
$$
\Vert(\phi-(1-\delta))^+\Vert_{L^\infty(\Omega\times(0,\widetilde{\tau}))}=0.
$$
We now repeat exactly the same argument for the case $(\phi-(-1+\delta))^-$ (using $\phi_n(t)=(\phi(t)+k_n)^-$), to get in the end that there exist $\delta>0,\widetilde{\tau}>0$ such that
\begin{align}
-1+\delta\leq \phi(x,t)\leq 1-\delta,\quad\text{a.e. in }\Omega\times (0,\widetilde{\tau}).
\label{endf}
\end{align}
Notice that $\widetilde{\tau}$ can be explicitly computed as a function of the parameters of the problem and the initial datum (see \eqref{tau2}). The proof is then concluded, recalling Remark \ref{proof} with $T_1=\widetilde{\tau}$.
\section{Some consequences of the Strict Separation Property}
\label{consequences}
In this Section we collect some results which are straightforward consequences of the strict separation property proven in Theorem \ref{sep}.
\subsection{Regularization in finite time}
 First we show that the weak solution given by Theorem \ref{known} actually regularizes more. Indeed, we have a first immediate consequence:
\begin{corollary}
	Under the same assumptions of Theorem \ref{sep}, for any $\tau>0$, there exists a constant $C=C(\tau)>0$ such that  
	\begin{align*}
	\Vert F^\prime(\phi(t))\Vert_{L^\infty(\Omega)}+\Vert\mu(t)\Vert_{L^\infty(\Omega)}\leq C,\qquad \forall t\geq \tau.
	\end{align*}
 \label{Linfty}
\end{corollary}
\begin{proof}
	The proof is immediate, since by the strict separation property we deduce $\Vert F^\prime(\phi(t))\Vert_{L^\infty(\Omega)}\leq C$ for any $t\geq \tau$ and then by comparison we get the $L^\infty$- control on $\mu$.
\end{proof}
Furthermore, we can also obtain the H\"{o}lder regularity of the weak solutions:
\begin{corollary}
	Under the same assumptions of Theorem \ref{sep}, for any $\tau>0$, there exists $C=C(\tau)>0$ and $\kappa=\kappa(\tau, \delta)\in (0,1)$ such that
	\begin{align}
	&\vert \phi(x_1,t_1)-\phi(x_2,t_2)\vert\leq C\left(\vert x_1-x_2\vert^\kappa+\vert t_1-t_2\vert^{\frac{\kappa}{2}}\right),\label{holder}\\&
	\vert \mu(x_1,t_1)-\mu(x_2,t_2)\vert\leq C\left(\vert x_1-x_2\vert^\kappa+\vert t_1-t_2\vert^{\frac{\kappa}{2}}\right)
	\label{holder2}
	\end{align}
	for all $(x_1,t_1),(x_2,t_2)\in \Omega_t$, where $\Omega_t=[t,t+1]\times\overline{\Omega}$ and $t\geq \tau$.
\end{corollary}
\begin{proof}
	We can argue as in \cite[Lemma 2.11]{GalGrasselli}. In particular, we rewrite the system \eqref{nonloc3d} in the following form 
	$$
	\partial_t\phi=\text{div}(a(x,\phi,\nabla\phi)),\quad (a(x,\phi,\nabla\phi)\cdot \textbf{n})_{\vert\partial\Omega}=0,
	$$
	with 
	$$
	a(x,\phi,\nabla\phi):=\Fsec(\phi)\nabla\phi-\nabla J\ast\phi.
	$$
	Since by assumption \ref{h1} we have $J\in W^{1,1}_{loc}(\R^3)$, $\Fsec(s)\geq \alpha$ for any $s\in(-1,1)$ by \ref{h3} \ry{and $\Vert \nabla J\ast \phi\Vert_{L^\infty(\Omega)}\leq \Vert \nabla J\Vert_{{L^1(B_r)}}$ by \eqref{J}, by Young's inequality we get
	\begin{align*}
	a(x,\phi,\nabla\phi)\cdot\nabla\phi&= \Fsec(\phi)\vert\nabla\phi\vert^2-(\nabla J\ast\phi)\cdot \nabla\phi\\&\geq \alpha\vert \nabla\phi\vert^2-\Vert \nabla J\Vert_{{L^1(B_r)}}\vert\nabla\phi\vert\\&
 \geq \frac{\alpha}{2}\vert \nabla\phi\vert^2-\frac{1}{2\alpha}\Vert \nabla J\Vert^2_{{L^1(B_r)}},
 \end{align*}
 and, similarly, by Corollary \ref{Linfty} and \eqref{J},
 $$
 \vert a(x,\phi,\nabla\phi)\vert\leq \Vert \Fsec(\phi)\Vert_{L^\infty(\Omega)}\vert \nabla\phi\vert+ \Vert \nabla J\ast \phi\Vert_{L^\infty(\Omega)}\leq C_1\vert \nabla\phi\vert+\Vert \nabla J\Vert_{{L^1(B_r)}},
	$$
	for some positive constant $C_1$ depending on $\tau, \delta$.} Therefore we infer the desired estimate \eqref{holder} applying \cite[Corollary 4.2]{Dung}. Then, by the regularity of $F$, we immediately deduce the same result for $\mu$, concluding the proof.
\end{proof}
In order to obtain higher-order spatial regularity for the phase variable $\phi$, we need to strengthen the assumptions on the interaction kernel $J$. In particular, we assume
\begin{enumerate}[label=(\subscript{H}{4})]
	\item \rx{Either} $J \in W^{2,1}(\mathcal{B}_{\rx{R}})$, where $\mathcal{B}_\rx{R}:=\{x\in \R^3: \vert x\vert<\rx{R}\}$, with $\rx{R}\sim\text{diam(}\Omega\text{)}$ such that $\overline{\Omega}\subset \mathcal{B}_\rx{R}$ \extra{and $x-\Omega\subset B_R$ for any $x\in \Omega$,} \rx{or $J$ is admissible in the sense of \cite[Definition 1]{Bedrossian}}.
	\label{h5}
\end{enumerate}
\begin{remark}
	As noticed in \cite[Remark 5.9]{GGG0}, we observe that Newtonian and \rx{second-order} Bessel potentials satisfy assumption \ref{h5}, \rx{namely they are admissible in the sense of \cite[Definition 1]{Bedrossian}}.
\end{remark}
Therefore we have
\begin{lemma}
		Under the same assumptions of Theorem \ref{sep}, assuming also that $J$ satisfies \ref{h5} and $F\in C^3(-1,1)$, for any $\tau>0$ there exists $C=C(\tau)>0$ such that 
		\begin{align}
	\Vert\phi\Vert_{L^{\frac{4}{3}}(t,t+1;H^2(\Omega))}\leq C, \quad \forall t\geq \tau.
		\label{H2}
		\end{align}
\end{lemma}
\begin{proof}
	We first observe that, since we can apply Theorem \ref{known}, by \eqref{Fprime}-\eqref{Lp} we deduce that 
	\begin{align}
		\Vert\nabla\phi\Vert_{L^\frac{8}{3}(t,t+1;L^4(\Omega))}+\Vert\mu\Vert_{L^2(t,t+1;V_2)}\leq C(\tau),\qquad \forall t\geq \tau,
	\label{maxcontrol}
	\end{align}
	for some positive constant $C(\tau)$. Then, as in \cite[Theorem 5]{Frigeri}, we proceed formally (these computations could be justified in a suitable approximating scheme, see, e.g., \cite[Theorem 5, Step 3]{Frigeri})  defining  $\partial_{ij}^2\phi:=\dfrac{\partial^2\phi}{\partial x_i\partial x_j}$, for $i,j=1,2,3$. We now apply $\partial_{ij}^2$ to the equation for the chemical potential $\mu$ and integrate on $\Omega$, to infer 
	\begin{align*}
	\int_\Omega\partial_{ij}^2\mu \partial_{ij}^2\phi dx=\int_\Omega \Fsec(\phi)(\partial_{ij}^2\phi)^2dx-\int_\Omega \partial_i(\partial_j J\ast \phi)\partial_{ij}^2\phi dx+\int_\Omega F^{\prime\prime\prime}(\phi)\partial_i\phi\partial_j\phi \partial_{ij}^2\phi,\qquad i,j=1,2,3.
	\end{align*}
	We now recall \rx{assumption \ref{h5}, so that by \cite[Lemma 2]{Bedrossian}:
	$$
	\Vert\partial_i(\partial_j J\ast \phi)\Vert_{L^2(\Omega)}\leq C\Vert\phi\Vert_{L^2(\Omega)}\leq C(\tau).
	$$}
	Therefore, by Cauchy-Schwartz and Young's inequalities, we infer, recalling that $\Fsec(s)\geq \alpha$ for any $s\in(-1,1)$, and exploiting the separation property of Theorem \ref{sep}, for any $t\geq \tau$,
	\begin{align*}
	\frac{\alpha}{2}\Vert\partial_{ij}^2\phi\Vert^2\leq C(1+\Vert\partial_{ij}^2\mu\Vert^2+\int_\Omega\vert\partial_i\phi\vert^2\vert\partial_j\phi\vert^2dx) \leq C\left(1+\Vert\mu\Vert^2_{H^2(\Omega)}+\Vert\nabla\phi\Vert^4_{L^4(\Omega)}\right),\quad i,j=1,2,3,
	\end{align*}
	which implies \eqref{H2}, thanks to \eqref{maxcontrol}.
\end{proof}
\subsection{Convergence to equilibrium}
We conclude the results of our paper by showing that the strict separation property is essential to study the longtime behavior of the single trajectory. In particular, we can follow \cite[Sec. 6.2]{DellaPorta}: for the sake of completeness we give here a sketch of the proofs. We employ the typical strategy based on the Lyapunov property of the associated system (see \eqref{dissipative}) and the well known \L ojasiewicz-Simon inequality. Let us consider the set of admissible initial data 
\begin{align*}
\mathcal{H}_m:=\left\{\phi\in L^\infty(\Omega): \Vert\phi\Vert_{L^\infty(\Omega)}\leq 1,\quad \vert\overline{\phi}\vert= m \right\}, 
\end{align*} 
with $m\in[0,1)$, and fix an initial datum $\phi_0\in \mathcal{H}_m$. Let then $\phi$ be the unique weak global-in-time solution departing from $\phi_0$, whose existence and uniqueness is ensured by Theorem \ref{known}. We introduce the $\omega$-limit set associated to $\phi_0$, i.e.,
\begin{align*}
\omega(\phi_0)=\{\widetilde{\phi}\in \mathcal{H}_m:\exists t_n\to \infty \text{ such that }\phi(t_n)\to \widetilde{\phi}\text{ in }H\}.
\end{align*}
By \eqref{H1}, $\phi$ is uniformly bounded in $V$, which is compactly embedded in $H$. Therefore, by standard results related to the intersection of non-empty, compact (in $H$), connected and nested sets, we infer that $\omega(\phi_0)$ is non-empty, compact and connected in $\mathcal{H}_m$. We now characterize the set $\omega(\phi_0)$, showing that it is composed by equilibrium points (i.e., stationary solutions) associated to \eqref{nonloc3d}, which are defined as 
\begin{definition}
	$\phi_\infty$ is an equilibrium point to problem \eqref{nonloc3d} if $\phi_\infty\in \mathcal{H}_m\cap V$ satisfies the stationary nonlocal Cahn-Hilliard equation 
	\begin{align}
		F^\prime(\phi_\infty)-J\ast \phi_\infty=\mu_\infty,\quad \text{in }\Omega,
		\label{conv1}
	\end{align}
	where $\mu_\infty\in \R$ is a real constant.
\end{definition}
As noticed also in \cite{DellaPorta}, the existence of a (not necessarily unique, see, e.g., \cite{Bates}) solution to \eqref{conv1} can be proven by means of a fixed point argument. Moreover, as shown in \cite[Lemma 6.1]{DellaPorta}, any $\phi_\infty\in V\cap \mathcal{H}_m$ satisfying \eqref{conv1} is strictly separated from the pure phases, i.e., there exists $\delta>0$ such that
$$
\Vert\phi_\infty\Vert_{L^\infty(\Omega)}\leq 1-\delta.
$$ 
If we now introduce the set of all the stationary points of the nonlocal Cahn-Hilliard equation: 
$$
\mathcal{S}:=\left\{\phi_\infty\in \mathcal{H}_m\cap V: \phi_\infty\text{ satisfies }\eqref{conv1}\right\},
$$
we can easily prove that $\omega(\phi_0)\subset \mathcal{S}$. Indeed, let us consider a sequence $t_n\to \infty$ such that $\phi(t_n)\to \widetilde{\phi}$ in $H$, $\widetilde{\phi}\in \omega(\phi_0)$. We then define the sequence of trajectories $\phi_n(t):=\phi(t+t_n)$ and $\mu_n(t):=\mu(t+t_n)$. Thanks to \eqref{H1}, we get, up to a non-relabeled subsequence, that 
$
\phi_n\overset{\ast}{\rightharpoonup} \phi^*\text{ in }L^\infty(0,\infty;V).
$
Passing to the limit in the equations for $\phi_n$, exploiting the results of Theorem \ref{known}, we infer that also $\phi^*$ satisfies \eqref{phi}-\eqref{mu} (we denote the corresponding chemical potential by $\mu^*$), with initial datum $\phi^*(0)=\widetilde{\phi}$. This last consideration follows from the fact that $\phi_n(0)=\phi(t_n)\to \widetilde{\phi}$ strongly in $H$. Moreover, we clearly have $\lim_{n\to \infty}\mathcal{E}(\phi_n(t))=\mathcal{E}(\phi^*(t))$ for all $t\geq 0$. By the energy identity \eqref{dissipative}, we infer that the energy $\mathcal{E}(\phi(\cdot))$ is nonincreasing in time, thus there exists $E_\infty$ such that $\lim_{t\to \infty}\mathcal{E}(\phi(t))=E_\infty$. This means that this convergence also holds for the subsequence $\{t+t_n\}_n$, thus 
$$\mathcal{E}(\phi^*(t))=\lim_{n\to \infty}\mathcal{E}(\phi_n(t))=\lim_{n\to \infty}\mathcal{E}(\phi(t+t_n))=E_\infty,$$
entailing that $\mathcal{E}(\phi^*(\cdot))$ is constant in time. Passing then to the limit in \eqref{dissipative}, which is valid for each $\phi_n$, we obtain 
$$E_\infty+\int_s^{t}\Vert\nabla\mu^*(\tau)\Vert^2d\tau\leq E_\infty,\quad \forall 0\leq s\leq t<\infty,$$
implying $\nabla\mu^*=0$ almost everywhere in $\Omega$, and thus, by comparison in \eqref{phi}, also $\partial_t\phi^*=0$ almost everywhere in $\Omega$, for almost every $t\geq 0$. Therefore, we infer that $\phi^*$ is constant in time, namely $\phi^*(t)=\widetilde{\phi}$ for all $t\geq 0$. \extra{Thus} $\widetilde{\phi}$ satisfies \eqref{conv1} with some constant $\mu_\infty\in \R$, and then $\widetilde{\phi}\in \mathcal{S}$, implying, being $\widetilde{\phi}\in \omega(\phi_0)$ arbitrary, $\omega(\phi_0)\subset \mathcal{S}$.
Notice that in this way we have shown that, for any $\phi_\infty\in \omega(\phi_0)$, 
\begin{align}
\mathcal{E}(\phi_\infty)=E_\infty=\lim_{s\to \infty}\mathcal{E}(\phi(s))=\inf_{s\geq 0}\mathcal{E}(\phi(s))\leq \mathcal{E}(\phi(t)), \qquad \forall t\geq 0.
\label{nonincr}
\end{align} 
We can then conclude by showing that $\omega(\phi_0)$ is a singleton. For the sake of clarity we present here the main tool, which is the \L ojasiewicz-Simon inequality (see, e.g., \cite[Proposition 6.2]{DellaPorta} or \cite{Loja}):
\begin{proposition}
	Let $P_0: H\to H_0$ be the projector operator. Assume that $F$ satisfies \ref{h3} and is real analytic in $(-1,1)$, $\phi\in V\cap L^\infty(\Omega)$ is such that $-1+\gamma\leq \phi(x)\leq 1-\gamma$, for any $x\in\overline{\Omega}$, for some $\gamma\in(0,1)$ and $\phi_\infty\in \mathcal{S}$. Then there exists $\theta\in \left(0,\frac{1}{2}\right)$, $\eta>0$ and a positive constant $C$ such that
	\begin{align}
	\vert\mathcal{E}(\phi)-\mathcal{E}(\phi_\infty)\vert^{1-\theta}\leq C\Vert P_0(F^\prime(\phi)-J\ast\phi)\Vert_{*},
	\label{ener}
	\end{align}
	whenever $\Vert\phi-\phi_\infty\Vert\leq \eta$.
	\label{Lojaw}
\end{proposition}
\begin{remark}
	We observe that the logarithmic potential $F$ is indeed real analytic in $(-1,1)$, thus the assumption of the foregoing proposition is satisfied.
\end{remark}
We have the following 
\begin{theorem}
	Under the same assumptions as in Theorem \ref{sep}, suppose additionally that $F$ is real analytic in $(-1,1)$. Then the weak solution $\phi$, departing from the initial datum $\phi_0\in \mathcal{H}_m$ converges to a single equilibrium point $\phi_\infty$ (depending on $\phi_0$) and $\omega(\phi_0)=\{\phi_\infty\}$. In particular we have 
	\begin{align}
	\lim_{t\to \infty}\Vert\phi(t)-\phi_\infty\Vert=0.
		\label{equil}
	\end{align} 
 \end{theorem}
\begin{proof}
	Thanks to \eqref{nonincr}, we infer that $\mathcal{E}(\phi(t))\geq \mathcal{E}(\phi_\infty)$, $\mathcal{E}(\phi(t))\to \mathcal{E}(\phi_\infty)$, as $t\to \infty$, for any $\phi_\infty\in \omega(\phi_0)$. Without loss of generality we can assume $\mathcal{E}(\phi(t))> \mathcal{E}(\phi_\infty)$ for all $t\geq 0$. Indeed, if there exists $\overline{t}>0$ such that $\E(\phi(\overline{t}))=\E(\phi_\infty)$, then clearly $\phi(t)=\phi(\overline{t})$ for any $t\geq \overline{t}$ and the claim follows, since then $\phi(t)=\phi_\infty$ for any $t\geq \overline{t}$. Let us now fix $\theta\in \left(0,\frac{1}{2}\right)$ and $\eta>0$ given in Proposition \ref{Lojaw}, where we have chosen $\gamma$ equal to the value of $\delta$ given in Theorem \ref{sep}.   By a contradiction argument as in \cite[Theorem 4]{Kreici} it is possible to show that there exists $t^*>0$ such that $\Vert\phi(t)-\phi_\infty\Vert\leq \eta$, for all $t\geq t^*$. Therefore, since the solution $\phi$ enjoys the separation property (by Theorem \ref{sep}) and thanks to the choice of $\gamma$, by Proposition \ref{Lojaw} we get, for any $t\geq t^*$,
	\begin{align*}
	\left(\E(\phi)-\E(\phi_\infty)\right)^{1-\theta}\leq\Vert P_0(F^\prime(\phi)-J\ast\phi)\Vert_{*}\leq C\Vert P_0\mu\Vert\leq \hat{C}\Vert\nabla\mu\Vert,
	\end{align*} 
where ${\hat{C}}>0$ depends on $C$ and on the Poincar\'{e}-Wirtinger's constant.
Therefore, by means of the energy identity \eqref{dissipative}, we deduce, for any $t\geq t^*$, 
\begin{align*}
&-\dfrac{d}{dt}\left(\E(\phi)-\E(\phi_\infty)\right)^\theta=-\theta\left(\E(\phi)-\E(\phi_\infty)\right)^{\theta-1}\frac{d}{dt}\mathcal{E}(\phi)\geq \dfrac{\theta\Vert\nabla\mu\Vert^2}{\hat{C}\Vert\nabla\mu\Vert}\geq \tilde{C}\Vert\nabla\mu\Vert,
\end{align*}
where $\tilde{C}>0$ is a positive constant independent of $t$. An integration over $(t^*,+\infty)$, for $t^*$ sufficiently large, implies that $\nabla\mu\in L^1(t^*,\infty;\textbf{H})$. By comparison, we deduce that also $\partial_t\phi\in L^1(t^*,\infty;V^\prime)$, so that  $$
\phi(t)=\phi(t^*)+ \int_{t^*}^t\partial_t\phi(\tau)d\tau\overset{t\to +\infty}{\longrightarrow}\widetilde{\phi}\quad \text{ in }V^\prime,
$$ 
for some $\widetilde{\phi}\in V^\prime$. Then we infer that $\phi(t)$ converges in $V^\prime$ as $t\to \infty$. By uniqueness of the limit in $V^\prime$, we can then conclude that $\omega(\phi_0)$ is a singleton, i.e., $\omega(\phi_0)=\{\widetilde{\phi}\}$. From now on we will denote $\widetilde{\phi}$ by $\phi_\infty$, since any $\phi_\infty\in \omega(\phi_0)$ coincides with $\widetilde{\phi}$. Thanks to \eqref{H1}, we then get \eqref{equil} by interpolation:
\begin{align*}
\Vert \phi(t)-\phi_\infty\Vert\leq C\Vert \phi(t)-\phi_\infty\Vert_V^{1/2} \Vert \phi(t)-\phi_\infty\Vert_{V^\prime}^{1/2}\leq C\Vert \phi(t)-\phi_\infty\Vert_{V^\prime}^{1/2}\overset{t\to +\infty}{\longrightarrow}0,
\end{align*}
concluding the proof.
\end{proof}
\subsection*{Acknowledgments} 
\extra{ The author thanks the anonymous referees for their appropriate comments
and useful remarks.}
The author is \extra{also} grateful to Andrea Giorgini and Maurizio Grasselli for several helpful comments on a preliminary version of this article and to Giorgio Meretti for the fruitful discussion about the Lemma on the Poincar\'{e}-type inequality. 
Moreover, the author has been partially funded by MIUR-PRIN Grant 2020F3NCPX \ \textquotedblleft Mathematics for Industry 4.0 (Math4I4)".

\end{document}